\documentclass[a4paper]{article}

\usepackage{graphicx}
\usepackage{amsmath}
\usepackage{amsfonts}
\usepackage{amssymb}
\newcommand{\be}{\begin{equation}}

\newcommand{\ee}{\end{equation}}

\newcommand{\bma}{\begin{pmatrix}}
\newcommand{\ema}{\end{pmatrix}}

\newcommand{\al}{\alpha}

\newcommand{\bet}{\beta}

\newcommand{\ga}{\gamma}

\newcommand{\Om}{\Omega}

\newcommand{\om}{\omega}

\newcommand{\G}{\Gamma}

\newcommand{\si}{\sigma}

\newcommand{\dz}{\wedge}

\newcommand{\ba}{\begin{array}}

\newcommand{\ea}{\end{array}}

\newcommand{\beq}{\begin{eqnarray}}

\newcommand{\eeq}{\end{eqnarray}}

\newtheorem{lm}{Lemma}

\newtheorem{thee}{Theorem}

\newtheorem{proo}{Proposition}

\newtheorem{co}{Corollary}

\newtheorem{rem}{Remark}

\newtheorem{deff}{Definition}

\newcommand{\bd}{\begin{deff}}

\newcommand{\ed}{\end{deff}}

\newcommand{\bl}{\begin{lm}}

\newcommand{\el}{\end{lm}}

\newcommand{\bp}{\begin{proo}}

\newcommand{\ep}{\end{proo}}

\newcommand{\bt}{\begin{thee}}

\newcommand{\et}{\end{thee}}

\newcommand{\bc}{\begin{co}}

\newcommand{\ec}{\end{co}}

\newcommand{\brm}{\begin{rem}}

\newcommand{\erm}{\end{rem}}

\newcommand{\der}{{\rm d}}

\begin{document}
\thispagestyle{empty}
\title{{Differential equations and conformal structures} 
\vskip 1.truecm
\author{
Pawe\l~ Nurowski\\
Instytut Fizyki Teoretycznej\\
Uniwersytet Warszawski\\
ul. Hoza 69, Warszawa\\
Poland\\
nurowski@fuw.edu.pl}
}
\maketitle

\begin{abstract}
We provide five examples of conformal geometries which are naturally
associated with ordinary differential equations (ODEs). The first example
describes a one-to-one correspondence between the Wuenschmann class of
3rd order ODEs considered modulo contact transformations of variables
and (local) 3-dimensional conformal Lorentzian geometries. The second
example shows that every point equivalent class of 3rd order ODEs
satisfying the Wuenschmann and the Cartan conditions define a
3-dimensional Lorentzian Einstein-Weyl geometry. The third example
associates to each point equivalence class of 3rd order ODEs a
6-dimensional conformal geometry of neutral signature. The fourth
example exhibits the one-to-one correspondence between point
equivalent classes of 2nd order ODEs and 4-dimensional conformal
Fefferman-like metrics of neutral signature. The fifth example shows
the correspondence between undetermined ODEs of the Monge type and
conformal geometries of signature $(3,2)$. The Cartan normal
conformal connection for these geometries 
is reducible to the Cartan connection with values
in the Lie algebra of the noncompact form of the exceptional group
$G_2$. All the examples are deeply rooted in Elie Cartan's works on
exterior differential systems. 
\end{abstract}

\newpage

\noindent

\rm
\section{Introduction}
\noindent
One aspect of the Null Surface Formulation of General Relativity (NSF)
of Fritelli, Kozameh and Newman \cite{NSF} is to encode the 
conformal geometry of space-time in the geometry of a certain pair of 
partial differential equations (PDEs) on the plane. Although this pair of 
differential equations appears in NSF quite naturally, the question
arises as to whether it is an accident or it is a feature of a deeper
link between differential equations and conformal structures. A closer
look at this question shows that the phenomenon observed in NSF is
only a tip of an iceberg, and that there is an abundance of examples
in which the geometry of differential equations can be related to the
conformal geometry in various dimensions. The main aim of this paper
is to describe these examples and to point out that all of them have
their roots in Elie Cartan's works on differential systems.\\

\noindent
The oldest and the simplest of these examples is due to Karl 
Wuenschmann. It is contained in his PhD dissertation \cite{Wun}
defended at the University of Greifswald in 1905. His result is quoted by 
Elie Cartan in a footnote of Ref.  
\cite{CartanSpanish}. According to Cartan Wuenschmann
observed that certain classes of 3rd order ordinary differential
equations (ODEs) define, in a natural way, a conformal Lorentzian
metric on the 3-dimensional spaces of their solutions. Chern in 
\cite{chern} interpreted the result of Wuenschmann in terms of a
Cartan normal conformal connection \cite{Kobayashi} with values in the Lie algebra 
$\underline{\bf so}(3,2)$. Recently, Newman and collaborators 
\cite{Newman1} proved that every 3-dimensional Lorentzian conformal 
geometry originates from a 3rd order ODE from the Wuenschmann class.\\

\noindent
Although, due to Cartan, we have the precise coordinates of Wuenschmann
thesis we were unable to get it from the University of
Greifswald. Thus, we do not know how Wuenschmann obtained his result. 
In a joint paper \cite{wunnewman} with Fritelli and Newman, we derived it by searching for 3rd ODEs for
which it was possible to define a null separation between the solutions. 
We believe, that this derivation is very close to the Wuenschmann
one. In the present paper, in Section 2, we give yet another
derivation of Wuenschmann's result. This presentation is 
closely related to the description given in Cartan's footnote. 
In particular, we specify under which differential
condition on $F=F(x,y,y',y'')$ the 3rd order ODE 
\be
y'''=F(x,y,y',y'')\label{ode33}
\ee
is in the Wuenschmann class (condition (\ref{wc})) and, using $F$ and
its derivatives, we give the explicit formula for 
the conformal Lorentzian 3-metric. We also calculate the conformal
invariants of these metrics, such as Cotton tensor, and relate them to
the contact invariants of the corresponding ODEs from the Wuenschmann
class. We end this section by providing nontrivial examples of ODEs
from the Wuenschmann class. \\

\noindent
Our next examples of conformal structures associated with differential
equations are motivated by Cartan's paper \cite{CartanSpanish}. In
this paper Cartan studies the geometry of an ODE (\ref{ode33}) given
modulo the point transformation of variables. He shows that if, in
addition to the Wuenschmann condition (\ref{wc}), the ODE satisfies
another point invariant condition (\ref{weylg}), then it defines a
3-dimensional Lorentzian Weyl geometry, i.e. the geometry defined by
a conformal class of Lorentzian 3-metrics $[g]$ and a 1-form $[\nu]$
given up to a gradient. This Weyl geometry turns out to satisfy
the Einstein-Weyl equations, which makes Cartan's observation important
in the integrable systems theory (see e.g. \cite{tod}).\\

\noindent
In Section 3 we formulate the equivalence problem for 3rd order
ODEs considered modulo point transformations and present its solution
(Theorem 3) due to Cartan \cite{CartanSpanish}. We 
interpret the result in terms of Cartan's connection $\om$ with values
in the Lie algebra of a 
group ${\bf CO}(1,2)\rtimes{\bf R}^3$ - the semidirect
product of the ${\bf SO}(1,2)$ group extended by the dilatations, and the
translation group in ${\bf R}^3$. In case of a generic 3rd order ODE 
(\ref{ode33}) this Cartan connection is defined on a principal  
${\bf SO}(1,2)$ fiber bundle $\cal P$ over a certain
{\it four} dimensional manifold but, if the equation satisfies the Wuenschmann
condition (\ref{wc}) and the Cartan condition (\ref{weylg}), it may be
interpreted as a Cartan connection on a principal ${\bf CO}(1,2)$ fiber
bundle over a {\it three} dimensional space identified with the
solution space of (\ref{ode33}). It is this special case which was
studied by Cartan. In
Section 3.1.1 we describe his result in the modern terminology. In
particular, we explicitly write down the formualae for the metric
$g_{ew}$ and the Weyl 1-form $\nu_{ew}$ in terms of function
$F=F(x,y,y',y'')$ defining the equation. We also prove that the
conditions (\ref{wc}) and (\ref{weylg}) are equivalent to the 
Einstein-Weyl condition for the Weyl geometry $[g_{ew},\nu_{ew}]$. The
result is summarized in Theorem 4. In two examples (Example 2 and
Example 3) we provide two nontrivial point equivalent classes of 3-rd
order ODEs which satisfy conditions (\ref{wc}) and
(\ref{weylg}). The class of equations of Example 2 is a generalization
of example (\ref{c2}) which was known to Tod \cite{tod}. Example 3
shows how to generate nontrivial $F=F(x,y,y',y'')$ satisfying
(\ref{wc}) and (\ref{weylg}) from particular solutions of 
reductions of the Einstein-Weyl geometries in 3-dimensions. Even very 
simple solutions, such as a solution $u=\sqrt{2x}$ of the dKP equation
(\ref{dkp1}), give rise to very nontrivial $F$s (see formula
(\ref{fdkp2})).\\

\noindent
In Section 3.1.2 we return to the generic case of an ODE (\ref{ode33})
given modulo point transformation and its Cartan connection $\om$ on 
the ${\bf SO}(1,2)$ fiber bundle $\cal P$. We show that in this
general case $\cal P$ is equipped with a special vector field whose
integral curves foliate $\cal P$. The 6-dimensional space of leaves of this
foliation is naturally equipped with a conformal metric
$[\tilde{\tilde{g}}]$ of signature
$(3,3)$. This 6-dimensional conformal structure encodes all the point
invariant information about the point equivalent class of ODEs
(\ref{ode33}). In particular, the Cartan (point) invariants of Theorem
3 and the curvature of Cartan's connection $\om$ can be equivalently 
described in terms of a Cartan normal conformal 
connection associated with the conformal class of metrics
$[\tilde{\tilde{g}}]$. This result, which was not mentioned by Cartan,
is summarized in Theorem 5; an explicit formula for this normal
conformal connection is given by (\ref{caln}).\\

\noindent
Section 4 deals with a geometry of a 2nd order ODE
\be
y''=Q(x,y,y')\label{ode22}
\ee 
considered modulo point transformations of variables. It provides a
next example of appearances of conformal geometry in the theory of
differential equations. This case was studied by us in a joint paper
with Sparling \cite{nurspar}. In this paper, exploiting an analogy
between 2nd order ODEs and 3-dimensional CR-structures, we were able
to associate a conformal 4-metric of signature (2,2) with each point
equivalence class of ODEs (\ref{ode22}). The construction of this metric,
described in \cite{nurspar}, was motivated by Fefferman's construction 
\cite{Fef} of Lorentzian metrics on a circle bundle over nondegenerate
3-dimensional CR-structures. Cartan, who formulated and solved the
equivalence problem for ODEs (\ref{ode22}) given modulo point
transformations in his famous paper \cite{Cartpc}, overlooked existence
of this metric. In Ref. \cite{nurspar}, we showed that the conformal
class of Fefferman-like (2,2) signature metrics associated with a
point equivalence class of ODEs (\ref{ode22}) encodes all the point
invariant information of such class. We summarize these results in
Theorem 6.\\

\noindent
In Section 5 which, in our opinion, is the highlight of the paper, we consider
the geometry of an undetermined equation
\be
z'=F(x,y,y',y'',z)\label{hilb11}
\ee    
for two real functions $y=y(x)$ and $z=z(x)$ of one variable. The
studies of equations of this type can be traced back to Gaspard Monge,
who knew that every solution to the equation of the form
\be
z'=F(x,y,y',z)\label{mon11}
\ee 
was expressible {\it without integrals} by means of an arbitrary
function of one variable and its derivatives. Hilbert 
\cite{hilbert}, on an example of equation 
\be
z'=(y'')^2,
\label{hilex}
\ee
showed that, in general, equations (\ref{hilb11}) do not have this
property. This result impressed Cartan, who previously
(Ref. \cite{5var}) considered equations
(\ref{hilb11}) as equations for Cauchy characteristics
of pairs of PDEs in the involution defined on the plane. Cartan solved
the equivalence problem for these PDEs which, implicitly, solves
an associated equivalence problem for ODEs (\ref{hilb11}). From
Cartan's solution of this equivalence problem it follows that among
equations (\ref{hilb11}) only those for which $F_{y''y''}=0$ have
general solutions which can be expressed {\it without integrals}.\\

\noindent
From the geometric point of view equations (\ref{hilb11}) for which
$F_{y''y''}\neq 0$ are much more interesting then those with
$F_{y''y''}=0$. It follows from Cartan's work \cite{5var} that
nonequivalent classes of equations (\ref{hilb11}) with $F_{y''y''}\neq
0$ are distinguished by means of a curvature of a certain 
Cartan connection. Surprisingly, this connection has 
values in the Lie algebra of the noncompact form $\tilde{G}_2$ of the {\it exceptional}
group $G_2$. The curvature of this connection is vanishing precisely
in the case of equations equivalent to the Hilbert example
(\ref{hilex}). This, in particular means that the symmetry group of a
very simple equation (\ref{hilex}) is isomorphic to
$\tilde{G}_2$. This fact, noticed with pride by Cartan in \cite{5var}, 
was perhaps the first geometric realization of this group predicted to exist by
Cartan and Engel in 1894.\\

\noindent
The main original part of Section 5 consists in an observation that
this Cartan connection can be understood as a reduction of a certain
Cartan normal {\it conformal} connection. This is associated with a conformal metric
$G_{(3,2)}$ of signature $(3,2)$ naturally defined by (\ref{hilb11})
on a 5-dimensional
space $J$ parametrized by the five independent variables
$(x,y,y',y'',z)$. It follows that all the invariant information about
the ODE (\ref{hilb11}) satisfying $F_{y''y''}\neq 0$ is encoded in the
conformal properties of the metric $G_{(3,2)}$. We introduced this
metric motivated by the Fefferman construction described in
\cite{nurspar}. Surprisingly, its existence, like the existence of 
Fefferman-like metrics described in Section 4, was overlooked by Cartan.\\

\noindent
Section 5 has three subsections. The first one makes precise the
notion of an equation having a general solution without integrals. It
also contains the proof of Monge's result on equation (\ref{mon11})
quoted above. The proof uses
Cartan's method of equivalence \cite{Olver} and aims to motivate the
definition of equivalence problem for equations (\ref{hilb11}). This
definition is given in Section 5.2 in terms of an equivalence of a 
system of three 1-forms (\ref{fhil}) on $J$. The beginning of Section 5.3 
reformulates 
Cartan's solution for the equivalence problem for
pairs of PDEs in involution on the plane adapting it to the
equivalence problem for ODEs (\ref{hilb11}) with $F_{y''y''}\neq
0$. This is summarized in Theorem 8. The interpretation of this
result in terms of Cartan's $\underline{\tilde{g}_2}$-valued
connection $\om_{\tilde{G}_2}$ 
is given by formula (\ref{ccg2}). The rest of this section
is devoted to the introduction and the discussion of a 5-dimensional
conformal $(3,2)$-signature metric whose Cartan normal conformal
connection is reducible to $\om_{\tilde{G}_2}$. This metric is defined
by formula (\ref{met32}) and is finally expressible entirely in terms
of the function $F=F(x,y,y',y'',z)$ and its derivatives in formula
(\ref{32met}). The main properties of this metric are summarized in
Theorem 9.\\

\noindent
As an application of this section, in Example 6, we
consider equations of the form 
$$
z'=F(y'').
$$
This generalizes (\ref{hilex}). We show that in this case there is
only one basic invariant of such equations. The metrics $G_{(3,2)}$ of
Example 6 turn out to be always
conformal to Einstein metrics. We characterize the Einstein scale for
them by means of a simple ODE. Finally, in case of a generic $F$, 
we show that the square of the Weyl tensor for metrics
$G_{(3,2)}$ can be interpreted in terms of a classical invariant of a
certain polynomial of the fourth order. This polynomial 
resembles very much the Weyl tensor 
polynomial known in the Newman-Penrose formalism \cite{NP}.

\section{Third order ODEs considered modulo contact transformations}
In 1905 Wuenschmann \cite{Wun} observed that the spaces of solutions of a
certain class of 3rd order ODEs are naturally equipped with conformal
Lorentzian geometries. His observation can be summarized as follows.\\

\noindent
Consider a 3rd order ordinary differential equation 
\be
y'''=F(x,y,y',y''),\label{ode3}
\ee
for a real function $y=y(x)$ of one variable. To simplify notation let 
$p=y'$ and $q=y''$. Now, consider the four-dimensional space $J^2$
parametrized by $(x,y,p,q)$. This space, the second jet space, is a
natural arena to study the geometry of equation (\ref{ode3}). In
particular, the total differential vector field 
$$
{\cal D}=\partial_x+p\partial_y+q\partial_p+F\partial_q
$$  
on $J^2$ yields the basic information about the solutions of
(\ref{ode3}). The integral curves of ${\cal D}$ foliate $J^2$ with
1-dimensional leaves. The leaf space 
${\cal S}$ of this foliation is 3-dimensional and can be identified
with the 3-dimensional space of solutions of
(\ref{ode3}). Following Chern \cite{chern} we equip $J^2$ with the following
bilinear form\footnote{Here and in the following we adapt the
  convention from General Relativity in which a symmetrized tensor
  product of two 1-forms $\al$ and $\beta$ is denoted by
  $\al\beta=\frac{1}{2}(\al\otimes\beta+\beta\otimes\al)$, e.g. $\al^2=\al\otimes\al$.} 
\be
\tilde{g}=2~[~\der y-p\der x~]~[~\der q-\frac{1}{3}F_q\der p+K\der
  y+(\frac{1}{3}q F_q-F-pK)\der x~]~-~[~\der p -q\der x~]^2.\label{metw}
\ee 
where
$$
K\equiv\frac{1}{6}{\cal D}F_q-\frac{1}{9}F_q^2-\frac{1}{2}F_p.
$$
Clearly, this form is degenerate. It has signature $(+,-,-,0)$ and its  
degenerate direction is tangent to the vector field ${\cal D}$. It is natural
to ask about the transformation properties of $\tilde{g}$ under the Lie
transport along the degenerate direction ${\cal D}$. It follows that $\tilde{g}$
transforms conformally under the Lie transport along ${\cal D}$ if and only
if the function $F=F(x,y,p,q)$ defining the ODE satisfies the
following nonlinear differential condition 
\be
A\equiv F_y~+~({\cal D}-\frac{2}{3}F_q)~K~=0.\label{wc}
\ee
This condition, {\it the Wuenschmann condition} defines {\it the
  Wuenschmann class} of 3rd order ODEs. Each equation from this class 
has a naturally defined conformal Lorentzian structure on the space of
its solutions. In our description, if $F$ satisfies
(\ref{wc}), this structure is obtained by projecting $\tilde{g}$ from $J^2$ to
the leaf space ${\cal S}$ of integral lines of ${\cal D}$. Since in such case $\tilde{g}$
transforms conformally along ${\cal D}$, it projects to the conformal
$(+,-,-)$ signature structure $[g]$ on ${\cal S}$. An interesting 
feature of the
above Wuenschmann construction is its invariance under the contact
transformations of the ODE (\ref{ode3}). More precisely, if equation
(\ref{ode3}) undergoes a transformation of variables of the form 
$$
x\to\bar{x}=\bar{x}(x,y,p),~~~y\to\bar{y}=\bar{y}(x,y,p),~~~p\to\bar{p}=\bar{p}(x,y,p),
$$
with
$$
\bar{y}_x-\bar{p}\bar{x}_x+p(\bar{y}_y-\bar{p}\bar{x}_y)=\bar{y}_p-\bar{p}\bar{x}_p=0
$$
then, if it is in the Wuenschmann class for the
function $F=F(x,y,p,q)$, it is also in the Wuenschmann class for the transformed
function $\bar{F}=\bar{F}(\bar{x},\bar{y},\bar{p},\bar{q})$. It
follows from the work of Chern \cite{chern} that the Wuenschmann
condition is the lowest order contact invariant condition one can
build out of $F$ and its partial derivatives. Moreover, every other
contact invariant of an equation from the Wuenschmann class
corresponds to a conformal invariant of the Lorentzian conformal
structure $[g]$. These conformal invariants are constructed by means
of the derivatives of the Cotton tensor $C$ of $[g]$. Assuming $A=0$ and
using the explicit form of the projection $[g]$ of $\tilde{g}$ we calculate
that the five independent components of $C$ are
\beq
&C_1=F_{qqqq},~~~~~~~~~~C_2=K_{qqq},~~~~~~~~~~C_3=L_{qq},~~~~~~~~~~C_4=N_q\nonumber\\
&~\nonumber\\
&C_5=-3 K_{qq} L + 3 K_q L_q - 3 K L_{qq} + 3 L_{qy} + 3 N_p + 
F_q N_q,\nonumber
\eeq
where
\beq
&L\equiv-\frac{1}{3}F_{qy}+\frac{1}{3}F_{qq}K-K_p-\frac{1}{3}F_q K_q\nonumber\\&~\nonumber\\
&N\equiv\frac{1}{3}F_{qq}L-\frac{2}{3}F_qL_q-2L_p+K K_{qq}-K_{qy}-\frac{1}{2}K_q^2.\nonumber
\eeq
It is worth noting that the vanishing of $C_1$ implies the vanishing
of all the $C_i$s, so that the conformal structure $[g]$ has vanishing
Cotton tensor iff $F_{qqqq}=0$. In such a case the corresponding
Wuenschmann class of equations (\ref{ode3}) is contact equivalent to
the equation $y'''=0$.\\

\noindent
Summing up we have the following theorem.\\

\noindent
\bt {\rm (Wuenschmann)}\\
Third order ODEs of the form
$$y'''=F(x,y,y',y'')$$
split onto two main contact nonequivalent classes, the Wuenschmann
class and its complement. There are contact nonequivalent equations
within the Wuenschmann class. Each representative of a contact
equivalence class of
equations satisfying Wuenschmann condition defines a conformal
Lorentzian structure on the space of its solutions. The conformally
equipped solution spaces of contact equivalent equations are conformally
related, so that each equivalence class of equations for which $A=0$
has a natural 3-dimensional conformal Lorentzian structure associated with it.  
\et

\noindent
There is a converse to this theorem.\\

\bt {\rm (Frittelli, Kozameh, Newman)}\\
Every 3-dimensional Lorentzian conformal structure $[g]$ defines a
contact equivalence class of third order ODEs satisfying the
Wuenschmann condition.
\et

\noindent
Below, without the proof which can be found in \cite{Newman1}, we
sketch how to pass from $[g]$ to the associated class of ODEs.\\

\noindent
Given a conformal family of Lorentzian metrics $[g]$ on a 3-dimensional
manifold $\cal M$ we start with a particular representative $g$ of $[g]$. This, in local
coordinates $\{x^i\},~ i=1,2,3$, can be written as $g=g_{ij}\der x^i\der
x^j$. Since the metric $g$ is Lorentzian it is meaningful to consider
the eikonal equation
$$
g^{ij}\frac{\partial S}{\partial x^i}\frac{\partial S}{\partial
  x^j}=0,
$$
for the real-valued function $S=S(x^i)$ on $\cal M$. This equation,
being homogeneous in $S$, has the {\it complete} solution $S=S(x^i;s)$
depending on a single parameter $s$. Now, treating $x^i$s, $i=1,2,3$
as constant parameters and $s$ as an independent variable we 
eliminate $x^1, x^2$ and $x^3$ by triple differentiation of equation
$S=S(x^i;s)$ with respect to $s$. As a result we get a relation of
the form 
\be
S'''=F(s,S,S',S''),\label{ode3eik}
\ee
which shows that $S=S(s)$ satisfies an ODE of the 3rd order. It
follows that this equation satisfies the Wuenschmann condition (\ref{wc}). It
also follows that if we start with another representative $\bar{g}$ of $[g]$
and find the corresponding complete solution
$\bar{S}=\bar{S}(\bar{x}^i;\bar{s})$ of the corresponding 
eikonal equation we get a 3rd order ODE for
$\bar{S}=\bar{S}(\bar{s})$  which is
related to (\ref{ode3eik}) by a contact transformation of variables
$s,S$ and $S'$.\\

\noindent
{\bf Example 1}\\
It can be easily checked that 
$$
F(x,y,p,q)=\alpha\frac{[q^2+(1-p^2)^2]^{\frac{3}{2}}}{[1-p^2]^{\frac{3}{2}}}-3\frac{p
  q^2}{1-p^2}-p(1-p^2)
$$
satisfies the Wuenschmann condition (\ref{wc}) for all values of the 
real parameter $\alpha$. Moreover, the 3rd order ODEs
$y'''=F(x,y,y',y'')$ corresponding to different values of $\alpha>0$
are contact nonequivalent. It follows that the conformal Lorentzian
structures $[g]$ associated with such $F$s have 4-dimensional group of
conformal symmetries, which correspond to four contact symmetries of
the associated 3rd order ODE \cite{Godlinskiphd}.  

\section{Third order ODEs considered 
modulo point transformations}

Cartan \cite{CartanSpanish} considered 3rd order ODEs modulo point
transformations of variables. These transformations are more
restrictive then the contact transformations. They merely mix the
independent and dependent variables 
\be
x\to\bar{x}=\bar{x}(x,y),~~~~~~~~y\to\bar{y}=\bar{y}(x,y) \label{point}
\ee
of the equation (\ref{ode3}). Cartan in \cite{CartanSpanish} found a
full set of invariants which determine whether two 3rd order ODEs are
transformable to each other by a point transformation of variables. He
used his {\it equivalence method}. This method starts with a system of four
1-forms
\beq
&\omega^1=\der y-p\der x\nonumber\\
&\omega^2=\der p-q\der x\label{cf}\\
&\omega^3=\der q-F(x,y,p,q)\nonumber\der x\\
&\omega^4=\der x\nonumber
\eeq 
which an ODE of the form (\ref{ode3}) defines on the second jet space
$J^2$. Under transformations (\ref{point}) of the ODE (\ref{ode3}) the
forms (\ref{cf}) transform by
\beq
&\omega^1\to \bar{\omega}^1=\alpha\omega^1\nonumber\\
&\omega^2\to\bar{\omega}^2=\beta(\omega^2+\gamma\omega^1)\label{tcf}\\
&\omega^3\to\bar{\omega}^3=\epsilon(\omega^3+\lambda\omega^2+\mu\omega^1)\nonumber\\
&\omega^4\to\bar{\omega}^4=\nu(\omega^4+\sigma\omega^1),\nonumber
\eeq 
where $\alpha,\beta,\gamma,\epsilon,\lambda,\mu,\nu,\sigma$ are
functions on $J^2$ such that $\alpha\beta\epsilon\nu\neq 0$. These
functions are determined by each particular choice of point 
transformation (\ref{point}). Instead of working with forms
$(\omega^i)$, $i=1,2,3,4,$ which are defined on $J^2$ only up to
transformations (\ref{tcf}), Cartan considers a manifold parametrized
by $(x,y,p,q,\alpha,\beta,\gamma,\epsilon,\lambda,\mu,\nu,\sigma)$ and
forms 
\beq
&\theta^1=\alpha\omega^1\nonumber\\
&\theta^2=\beta(\omega^2+\gamma\omega^1)\nonumber\\
&\theta^3=\epsilon(\omega^3+\lambda\omega^2+\mu\omega^1)\nonumber\\
&\theta^4=\nu(\omega^4+\sigma\omega^1),\nonumber
\eeq 
which are well defined there. Using his equivalence method he
constructs a 7-dimensional manifold ${\cal P}$ on which the four forms
$\theta^1,\theta^2,\theta^3,\theta^4$ supplemented by three other forms
$\Omega_1, \Omega_2$ and $\Omega_3$ constitute a rigid coframe. This
coframe encodes all the point invariant information about the ODE
(\ref{ode3}). More precisely, Cartan proves the following theorem.
\bt~\\
A 3rd order ODE (\ref{ode3}) considered modulo point transformations of variables
(\ref{point}) uniquely defines 
\begin{itemize}
\item a 7-dimensional manifold ${\cal P}$, 
\item seven 1-forms
$\theta^1,\theta^2,\theta^3,\theta^4,\Omega_1,\Omega_2,\Omega_3$ on
  ${\cal P}$
such that
$\theta^1\dz\theta^2\dz\theta^3\dz\theta^4\dz\Omega_1\dz\Omega_2\dz\Omega_3\neq
0$ and
\item functions $A,B,C,D,G,H,K,L,M,N$ on ${\cal P}$, 
\end{itemize}
which satisfy the following
differential conditions
\beq
&\der\theta^1=\Omega_1\dz\theta^1+\theta^4\dz\theta^2\nonumber\\
&~\nonumber\\
&\der\theta^2=\Omega_2\dz\theta^2+\Omega_3\dz\theta^1+\theta^4\dz\theta^3\nonumber\\&~\nonumber\\
&\der\theta^3=(2\Omega_2-\Omega_1)\dz\theta^3+\Omega_3\dz\theta^2+A\theta^4\theta^1\nonumber\\&~\nonumber\\
&\der\theta^4=(\Omega_1-\Omega_2)\dz\theta^4+B\theta^2\dz\theta^1+C\theta^3\dz\theta^1\label{syspoint}\\&~\nonumber\\
&\der\Omega_1=-\Omega_3\dz\theta^4+(H+D)\theta^1\dz\theta^2+(3K-2B)\theta^1\dz\theta^3+(G+L)\theta^1\dz\theta^4-C\theta^2\dz\theta^3\nonumber\\&~\nonumber\\
&\der\Omega_2=D\theta^1\dz\theta^2+2(K-B)\theta^1\dz\theta^3+G\theta^1\dz\theta^4-2C\theta^2\dz\theta^3\nonumber\\&~\nonumber\\
&\der\Omega_3=(\Omega_2-\Omega_1)\dz\Omega_3+M\theta^1\dz\theta^2+(D-H)\theta^1\dz\theta^3+N\theta^1\dz\theta^4+(2K-B)\theta^2\dz\theta^3+G\theta^2\dz\theta^4.\nonumber
\eeq
Two third order ODEs $y'''=F(x,y,y',y'')$ and $\bar{y}'''=\bar{F}(\bar{x},\bar{y},\bar{y}',\bar{y}'')$ are transformable to
each other by means of a point transformation (\ref{point}) if and
only if there exists a diffeomorphism $\phi:{\cal P}\to\bar{\cal P}$ of the
corresponding manifolds $\cal P$ and $\bar{\cal P}$ such that
$\phi^*(\bar{\theta}^i)=\theta^i,$ $i=1,2,3,4,$ and
$\phi^*(\bar{\Omega}_\mu)=\Omega_\mu$, $\mu=1,2,3$.
\et 

\subsection{Cartan connections associated with 3rd order ODEs considered
  modulo point transformations}
Among the equivalence classes of 3rd order ODEs described by Theorem 3
there is a particularly simple class corresponding to the vanishing of
all the functions $A,B,C,D,G,H,K,L,M,N$. In case of such ODEs the
corresponding forms
$(\theta^1,\theta^2,\theta^3,\theta^4,\Om_1,\Om_2,\Om_3)$ can be
considered a basis of left invariant forms on a Lie group which
naturally identifies with the space $\cal P$.
The structure constants of this group are
determined by equations (\ref{syspoint}) with all the functions
$A,B,C,D,G,H,K,L,M,N$ vanishing. This group turns out to be
locally isomorphic to ${\bf CO}(1,2)\rtimes{\bf R}^3$, the semidirect
product of the ${\bf SO}(1,2)$ group extended by the dilatations, and the
translation group in ${\bf R}^3$. In this sense, Theorem 3 can be
interpreted in terms of a ${\bf CO}(1,2)\rtimes{\bf R}^3$ Cartan
connection defined over the space $J^2$. Explicitly, the 1-form 
\be
\omega=\bma 
\Omega_2&0&0&0&0\\
\theta^1&\Omega_2-\Omega_1&-\theta^4&0&0\\
\theta^2&-\Omega_3&0&-\theta^4&0\\
\theta^3&0&-\Omega_3&\Omega_1-\Omega_2&0\\
0&\theta^3&-\theta^2&\theta^1&-\Omega_2
\ema,\label{conpoint}
\ee
which has values in the Lie algebra of ${\bf CO}(1,2)\rtimes{\bf
  R}^3$, defines a Cartan connection on $\cal P$. To see this it is
enough to observe that the system (\ref{syspoint}) guarantees that the
annihilator of forms $(\theta^1,\theta^2,\theta^3,\theta^4)$ is
integrable, so that $\cal P$ is fibered over the 4-dimensional space of leaves
tangent to this annihilator. This space of leaves naturally identifies
with $J^2$. Using equations (\ref{syspoint}) and calculating 
$$
R=\der\omega+\omega\dz\omega
$$
to be 
$$
R=\bma {\cal
  F}&0&0&0&0\\0&R^1_{~1}&R^1_{~2}&0&0\\0&R^2_{~1}&0&R^1_{~2}&0\\\Theta^3&0&R^2_{~1}&-R^1_{~1}&0\\0&\Theta^3&0&0&-{\cal F}
\ema,
$$
with
\beq
&{\cal F}=D\theta^1\dz\theta^2+2(K-B)\theta^1\dz\theta^3-2C\theta^2\dz\theta^3+G\theta^1\dz\theta^4\nonumber\\
&~\nonumber\\
&\Theta^3=-A\theta^1\dz\theta^4\nonumber\\
&~\nonumber\\
&R^1_{~1}=-H\theta^1\dz\theta^2-K\theta^1\dz\theta^3-C\theta^2\dz\theta^3-L\theta^1\dz\theta^4\nonumber\\
&~\nonumber\\
&R^2_{~1}=-M\theta^1\dz\theta^2+(H-D)\theta^1\dz\theta^3+(B-2K)\theta^2\dz\theta^3-N\theta^1\dz\theta^4-G\theta^2\dz\theta^4\nonumber\\
&~\nonumber\\
&R^1_{~2}=B\theta^1\dz\theta^2+C\theta^1\dz\theta^3,\nonumber
\eeq
we find that ${\bf SO}(1,2)\to{\cal P}\to J^2$ equipped with $\omega$
is a Cartan bundle with a ${\bf CO}(1,2)\rtimes{\bf R}^3$ connection
over $J^2$.\\

\noindent
In the next subsection we discuss under which conditions Theorem 3 can
be interpreted in terms of a Cartan connection over a certain 
{\it 3-dimensional} space, the space with which all the solution
spaces of point equivalent equations (\ref{ode3}) may be identified.

\subsubsection{A subclass defining Lorentzian Einstein-Weyl geometries on the
  solution space}
\noindent
First, the system (\ref{syspoint}) guarantees that, not only 
$\cal P$ is foliated by the 3-dimensional leaves discussed so far, but it is 
also foliated by {\it 4-dimensional} leaves. These are tangent to the
integrable distribution on which the forms 
$(\theta^1,\theta^2,\theta^3)$ vanish. Thus $\pi:{\cal P}\to {\cal M}$
can be considered a fibre
bundle over the {\it 3-dimensional} space $\cal M$ of leaves of this
foliation. A 4-dimensional group ${\bf CO}(1,2)$ acts naturally on the fibres
$\pi^{-1}({\cal M})$
of $\cal P$ equipping it with a structure of a ${\bf CO}(1,2)$ fibre bundle
over $\cal M$. Now, the form $\om$ defined by (\ref{conpoint}) can be
interpreted as a ${\bf CO}(1,2)\rtimes{\bf R}^3$ Cartan connection on 
${\bf CO}(1,2)\to{\cal P}\to {\cal M}$ iff in the curvature $R$ there
are only horizontal $\theta^1\dz\theta^2$, $\theta^1\dz\theta^3$ and $\theta^2\dz\theta^3$ terms. This
is only possible if
\be
a)~~~~~~~~~~~~~~~~~~~~~~~~~~~~~~~~~~~~~~~
A\equiv 0~~~~~~~~~~~~~~~~~~~~~~~~~~~~~~~~~~~~~~~~~~~~~~~\label{wc'}
\ee
and
$$
b)~~~~~~~~~~~~~~~~~~~~~~~~~~~~~~~~~~~~~~~
G\equiv 0.~~~~~~~~~~~~~~~~~~~~~~~~~~~~~~~~~~~~~~~~~~~~~~
$$
These are also sufficient conditions since, if they are satisfied, 
the functions $N$ and $L$ also vanish. Vanishing
of each of $A$ and $G$ is a point invariant property of the ODE
(\ref{ode3}). One can also consider these conditions independently of each
other. The vanishing of $A$ is precisely the Wuenschmann condition
(\ref{wc}) which, being contact invariant, is also a point invariant. If
the equation (\ref{ode3}) satisfies this condition it defines the
conformal metric $g$ on $\cal M$. This conformal Lorentzian structure
on $\cal M$ is the projection of the bilinear form
\be
\tilde{g}=2\theta^1\theta^3-(\theta^2)^2\label{bf}
\ee
from $\cal P$ to $\cal M$. Thus, similarly to the contact case, point
equivalent classes of equations (\ref{ode3}) satisfying the
Wuenschmann condition $A=0$ define a conformal structure on the space
$\cal M$. If, in addition condition $b)$ is satisfied 
then the pair $(\tilde{g},\tilde{\nu}=-2\Omega_2)$ projects to a well 
defined {\it Weyl geometry} $[(g_{ew},\nu_{ew})]$ on the
space $\cal M$. We recall
  that a {\it Weyl geometry} on a 3-dimensional manifold $\cal M$ is the
geometry associated with an equivalence class $[(g,\nu)]$ of pairs
$(g,\nu)$, in which $g$ is a Lorentzian metric, $\nu$ is a 1-form, and
two pairs  $(g,\nu)$ and  $(g',\nu')$ are in the
equivalence relation iff there exists a function $\phi$ on $\cal M$ such that
$g'={\rm e}^{-2\phi}g$ and $\nu'=\nu+2\der\phi$.\\ To see how the Weyl
geometry $[(g_{ew},\nu_{ew})]$ appears in the above context we first remark that the condition $G\equiv 0$, when written 
in terms of the function $F=F(x,y,p,q)$ defining the equation
(\ref{ode3}), is
\be
G\equiv 0\quad\quad\quad\quad\iff\quad\quad\quad\quad {\cal
  D}^2F_{qq}-{\cal D}F_{qp}+F_{qy}=0.\label{weylg}
\ee 
Then, identifying $\cal M$ with the quotient $J^2/{\cal D}$ and using 
the $(x,y,p,q)$ coordinates on $J^2$, we have 
$$
\tilde{g}=\beta^2[2\om^1(\om^3-\frac{1}{3}F_q\om^2+(\frac{1}{6}{\cal D}F_q-\frac{1}{2}F_p-\frac{1}{9}F_q^2)\om^1)-(\om^2)^2]
$$
$$
-\tilde{\nu}=2\Omega_2=2\der\log\beta+\frac{2}{3}(F_{qp}-{\cal D}F_{qq})\om^1+\frac{2}{3}F_{qq}\om^2+\frac{2}{3}F_q\om^4.
$$
The bilinear form $\tilde{g}$ is identical with (\ref{metw}), thus due to the
Wuenschmann condition $A=0$, it projects to a conformal structure $[g_{ew}]$ on $\cal
M$. Calculating the Lie derivative of $\tilde{\nu}$ with respect to
${\cal D}$ we
find that
$$
{\cal L}_{\cal D}\tilde{\nu}=\frac{2}{3}({\cal D}^2F_{qq}-{\cal D}F_{qp}+F_{qy})\om^1+\der (...).
$$
Thus, due to condition (\ref{weylg}), ${\cal L}_{\cal D}\tilde{\nu}$ is a total
differential. This means that $\tilde{\nu}$ projects to a class of
1-forms $[\nu_{ew}]$ on $\cal M$ which are given up to an addition of a
gradient.\\
It follows that the so defined Weyl geometry $[(g_{ew},\nu_{ew})]$ on $\cal M$ 
satisfies the Einstein-Weyl equations. To see this we first 
recall that a 3-dimensional Weyl geometry 
$[(g=g_{ij}\theta^i\theta^j,\nu)]$ defines a Weyl connection, which
is totally determined by the connection 1-forms
$\G^i_{~j}$ satisfying 
$$\der\theta^i+\G^i_{~j}\dz\theta^j=0,~~~~~
\der g_{ij}-\G_{ij}-\G_{ji}=-\nu g_{ij},\quad\quad\Gamma_{ij}=g_{ik}\Gamma^k_{~j}.$$
The Weyl geometry is said to be Einstein-Weyl iff the curvature 
$$\Om^i_{~j}=\frac{1}{2}
R^i_{~jkl}\theta^k\dz\theta^l=
\der\G^i_{~j}+\G^i_{~k}\dz\G^k_{~j}$$
of this connection satisfies 
\be R_{(ij)}-\frac{1}{3}Rg_{ij}=0,\label{eww}\ee
where $$R_{jk}=R^i_{~jik} \quad\quad R_{(ij)}=\frac{1}{2}(R_{ij}+R_{ji})
$$ and  $R=g^{ij}R_{ij}$, $g^{ik}g_{kj}=\delta^i_{~j}$.\\
It follows that in the case of Weyl geometry $[(g_{ew},\nu_{ew})]$ 
the Einstein-Weyl condition (\ref{eww}) reduces to the requirement
that the point invariant $M$ of the system (\ref{syspoint}) vanishes. 
To show that conditions $A=G=0$, which were needed to define
$[(g_{ew},\nu_{ew})]$, imply $M=0$ we apply the exterior derivative d
to the both sides of equations (\ref{syspoint}). Then from the
equation $\der^2\theta^3=0$ we deduce that $N=L=0$. Having this and 
insisting on $\der^2\Om_2=0$ we get that $D=2H$, which is only
possible if $M=0$.\\

\noindent
Summarizing we have the following theorem.

\bt {\rm (Cartan)}\\
A point equivalence class of 3rd order ODEs represented by an ODE
$$
y'''=F(x,y,y',y'')
$$
which satisfies Wuenschmann condition (\ref{wc}) and Cartan condition 
(\ref{weylg}) defines a Lorentzian Einstein-Weyl geometry 
$[(g_{ew},\nu_{ew})]$ on the
3-dimensional space $\cal M$. This space can be identified with the 
solution space of any of the ODEs from the equivalence
class.
\et

\noindent
It is a nontrivial task to find $F=F(x,y,p,q)$ which 
satisfies the Einstein-Weyl conditions (\ref{wc}), (\ref{weylg}). Cartan gave several examples of such $F$s 
(see \cite{tod} for a discussion of that issue). Here we present two
other ways of constructing them.\\

\noindent
{\bf Example 2}\\
It is relatively easy to find all point equivalence classes of 3rd
order ODEs which admit at least {\it four} infinitesimal 
point symmetries \cite{Godlinskiphd}. Among them there is a
1-parameter family of nonequivalent ODEs represented by 
\be
F=\frac{(\sqrt{a(2qy-p^2)}~)^3}{y^2},\label{c1}
\ee 
which corresponds to nonequivalent 
Einstein-Weyl geometries for each value of the real constant
$a$. This constant enumerates nonequivalent ODEs; its sign is
correlated with the sign of $(2qy-p^2)$, so that the expression under
the square root is positive. If $a\to\infty$ the equivalence class of
ODEs may be represented by 
\be
F=q^{3/2},\label{c2}
\ee
which also satisfies conditions (\ref{wc}), (\ref{weylg}). \\

\noindent
{\bf Example 3}\\
Since 3-dimensional Lorentzian Einstein-Weyl geometries are known to
be generated by solutions of various integrable systems, one can try to
use such solutions to associate with them point equivalence classes of
ODEs (\ref{ode3}). We illustrate this procedure on an
example of solutions to the dKP equation.\\

\noindent
The dKP equation for a real function $u=u(x,y,t)$ can be considered to
be the Froebenius condition 
\beq
&\der\bar{\om}^{1}\dz\bar{\om}^{1}\dz\bar{\om}^{4}=0\nonumber\\
&\der\bar{\om}^{4}\dz\bar{\om}^{1}\dz\bar{\om}^{4}=0\label{dkp}
\eeq
for the two Pfaffian forms 
\beq
&\bar{\om}^{1} = \der x +(u+v^2)\der t + v \der y\nonumber\\
&\bar{\om}^{4} = \der v - (u_y + u_x v) \der t  - u_x \der y\label{fdkp}
\eeq
in a 4-dimensional space parametrized by $(x,y,t,v)$. Indeed, by
substitution of (\ref{fdkp}) to (\ref{dkp}) we find that (\ref{dkp})
is equivalent to 
\be
u_{yy}=-(u_x)^2+u_{xt}-u u_{xx},\label{dkp1}
\ee
which is the dKP equation. Since every solution to (\ref{dkp1})
generates a 3-dimensional Lorentzian Einstein-Weyl geometry 
\cite{dunajski} it is 
reasonable to ask if there is a point equivalence class of 3rd order
ODEs associated with each such solution. It turns out that the answer
to this question is positive. Given a solution $u=u(x,y,t)$ of the 
dKP equation there is a point equivalence class of 3rd order ODEs, with
a representative in the form (\ref{ode3}), such that the four 1-forms
$(\bar{\om}^1,\bar{\om}^2,\bar{\om}^3,\bar{\om}^4)$ of (\ref{tcf}) encoding it have, in a convenient coordinate system
$(x,y,t,v)$ on
$J^2$, representatives $\bar{\om}^1$ and $\bar{\om}^4$ of (\ref{fdkp})
and $\bar{\om}^2$ and $\bar{\om}^3$ given by 
\beq
&\bar{\om}^2 = (-u u_{xx} - 2 u_{xy} v + u_{xx} v^2)\der t - u_{xx}\der
x - u_{xy} \der y\nonumber\\
&\bar{\om}^3 = (-u u_{xx}^2 - 4 u_{xy}^2 + 4 u_{xx}u_{xy} v - u_{xx}^2
v^2)\der t - u_{xx}^2\der x + u_{xx}(-2 u_{xy} + u_{xx}v)\der y.\nonumber
\eeq 
In particular, equations (\ref{dkp}) guarantee that there exists a
coordinate $X$ on $J^2$ such that in the class (\ref{tcf}) of forms
$\bar{\om}^4$ there is an exact form $\der X$. This defines a function
$X$, which in turn is
interpreted as the independent variable of the associated ODE. For
example, for a very simple solution
$$
u=\sqrt{2x}
$$
of the dKP equation we find that $$X=t+\frac{1}{2}v^2+\sqrt{2x},$$ which enables us
to find the associated class of 3rd order ODEs. This class may be 
represented by quite a nontrivial 
\be
F(x,y,p,q)=\frac{p q~(-12 + 3 p q - 8 \sqrt{1 - p q}) + 8(1 + \sqrt{1 - p
 q})}{p^3}.\label{fdkp2}
\ee
It can be checked by a direct substitution that such $F$ satisfies 
the Einstein-Weyl conditions (\ref{wc}), (\ref{weylg}).\\

\noindent
We close this section with a remark, that it is not clear whether all
3-dimensional Lorentzian Einstein-Weyl geometries have their associated
point equivalence classes of 3rd order ODEs. Our experience, based on 
the Cartan's equivalence method, suggests that it is very likely.    

\subsubsection{Conformal metric of signature $(3,3)$ associated
  with a point equivalence class of 3rd order ODEs}

\noindent
If an ODE (\ref{ode3}) does not satisfy the Wuenschmann condition
(\ref{wc'}), it is
impossible to define a conformal structure in 3 dimensions out of the
Cartan invariants (\ref{syspoint}). However, irrespectively of the
Wuenschmann condition (\ref{wc'})
being satisfied or not, with each point equivalence class of ODEs
(\ref{ode3}), we can associate a conformal metric of signature
$(+,+,+,-,-,-)$, whose conformal invariants encode all the point
invariant information about the corresponding class of ODEs. We
achieve this by using Sparling's procedure \cite{sparlingpriv} which, with 
`the Levi-Civita part'
$$
\G=(\G^i_{~j})=\bma 
\Omega_2-\Omega_1&-\theta^4&0\\
-\Omega_3&0&-\theta^4\\
0&-\Omega_3&\Omega_1-\Omega_2\\
\ema,
$$
of the Cartan connection (\ref{conpoint}) and with the bilinear form 
$\tilde{g}=g_{ij}\theta^i\theta^j$ of (\ref{bf}), associates a new 
bilinear form 
$$
\tilde{\tilde{g}}=\epsilon_{ijk}\theta^i\G^j_{~l}g^{lk}
$$
on $\cal P$. Here $$(g^{ij})=\bma 0&0&1\\0&-1&0\\1&0&0\\\ema$$ and $\epsilon_{ijk}$ is the standard
Levi-Civita symbol in ${\bf R}^3$ so that 
$$
\tilde{\tilde{g}}=2[~(\Om_1-\Om_2)\theta^2-\Om_3\theta^1+\theta^4\theta^3~].
$$
This bilinear form is degenerate on $\cal P$ and has $(+,+,+,-,-,-,0)$ 
signature. Denoting the basis of vector fields on $\cal P$ dual to the 1-forms
$(\theta^1,\theta^2,\theta^3,\theta^4,\Om_1,\Om_2,\Om_3)$ by
$(X_1,X_2,X_3,X_4,Y_1,Y_2,Y_3)$, we find that the degenerate direction of
$\tilde{\tilde{g}}$ is tangent to the vector field $Z=Y_1+Y_2$.\\ 

\noindent
It is
remarkable that, due to equations (\ref{syspoint}), the bilinear form
$\tilde{\tilde{g}}$ transforms conformally when Lie transported along
  $Z$. Explicitly, without any assumptions on the Cartan invariants
  $A,B,C,D,$ $G,H,K,L,M,N$, we have 
$${\cal L}_Z\tilde{\tilde{g}}=\tilde{\tilde{g}}.$$
Thus, the bilinear form $\tilde{\tilde{g}}$ naturally descends to a 
conformal metric $g_{\cal N}$ of neutral signature on 
the 6-dimensional space $\cal N$ of integral curves of the
vector field $Z$. This conformal metric yields all the point invariant
information about the corresponding point equivalent class of ODEs
(\ref{ode3}). In particular, the Cartan invariants
$A,B,C,D,G,H,K,L,M,N$ can be understood as curvature coefficients of
the Cartan normal conformal connection associated with $g_{\cal
  N}$. This Cartan connection can be represented by the following ${\bf
  so}(4,4)$-valued 1-form 
\be
\om_{\cal N}=\bma
\frac{1}{2}\Om_2&\frac{1}{4}(\Om_1 -
\Om_2)&-\frac{1}{4}\theta^4&\frac{1}{4}\Om_3&\tau_4&\tau_5&\frac{1}{2}\G^3_{~4}&
0\\&&&&&&&\\
\Om_1 - \Om_2&\frac{1}{2}\Om_2&\frac{1}{2}\theta^4&\G^1_{~3}& 0&-\G^2_{~4}&-\G^3_{~4}&\tau_4\\&&&&&&&\\
-\Om_3&\frac{1}{2} \Om_3&\frac{1}{2}\Om_1&\G^2_{~3}&\G^2_{~4}& 0& \G^2_{~6}&\tau_5\\&&&&&&&\\
\theta^4&\frac{1}{2}\theta^4& 0&-\frac{1}{2}\Om_1 + \Om_2&
\G^3_{~4}&-\G^2_{~6}& 0& \frac{1}{2}\G^3_{~4}\\&&&&&&&\\
\theta^2& 0&-\frac{1}{2}\theta^1&
\frac{1}{2}\theta^3&-\frac{1}{2}\Om_2& -\frac{1}{2}\Om_3&
-\frac{1}{2}\theta^4&\frac{1}{4}(\Om_1 - \Om_2)\\&&&&&&&\\
\theta^1& \frac{1}{2}\theta^1& 0& \frac{1}{2}\theta^2&
-\frac{1}{2}\theta^4& -\frac{1}{2}\Om_1& 0&-\frac{1}{4}\theta^4\\&&&&&&&\\
\theta^3& -\frac{1}{2}\theta^3&-\frac{1}{2}\theta^2&
0&-\G^1_{~3}& -\G^2_{~3}& \frac{1}{2}\Om_1 -
\Om_2&\frac{1}{4}\Om_3\\&&&&&&&\\
0&\theta^2&\theta^1&\theta^3&\Om_1 - \Om_2& -\Om_3& \theta^4&-\frac{1}{2}\Om_2
\ema,\label{caln}
\ee 
where 
\beq
&\tau_4=\frac{1}{12}[X_3(G) - 6 H]~\theta^1
  - \frac{1}{4}K~\theta^2 - \frac{1}{2}C~\theta^3\nonumber\\&\tau_5=\frac{1}{2}[-A C - 2
  X_2(L) - 2 M + X_4(D)]~\theta^1 +\frac{1}{12}[X_3(G) - 6 H + 6 D]~\theta^2 +
  \frac{1}{4}(-2 B + 3 K)~\theta^3 + \frac{1}{2}G~\theta^4\nonumber\\
&\G^1_{~3}=\frac{1}{2}\Om_3 +
\frac{1}{2}(G - L)~\theta^1\nonumber\\
&\G^2_{~3}= N~\theta^1 + \frac{1}{2}(G -
L)~\theta^2 + A~\theta^4\nonumber\\
&\G^2_{~4}= M~\theta^1 - H~\theta^2 + \frac{1}{2}(2 B - 3
K)~\theta^3 - \frac{1}{2}(G + L)~\theta^4\nonumber\\
&\G^2_{~6}=(-H +D )~\theta^1 +
\frac{1}{2}K~\theta^2 + C~\theta^3\nonumber\\
&\G^3_{~4}=\frac{1}{2}
(2 B - K)~\theta^1 +C~\theta^2\nonumber
\eeq
on $\cal P$.

\noindent
We remark that not all 6-dimensional split-signature conformal metrics
originate from a point equivalence 3rd order ODEs. To see this, we
calculate the curvature 
$$R_{\cal N}=\der\om_{\cal N}+\om_{\cal N}\dz\om_{\cal N}$$
of $\om_{\cal N}$ and observe\footnote{We omit writing down 
the explicit formulae for this curvature here.} that it has quite special form when
compared to the curvature of Cartan's normal conformal connection
associated with a generic $(+,+,+,-,-,-)$ signature metric. \\

\noindent
Summarizing, we have the following theorem.
\bt~\\
Each point equivalence class of 3rd order ODEs $$y'''=F(x,y,y',y'')$$
defines a conformal split-signature metric $g_{\cal N}$ on a 
6-dimensional manifold $\cal N$, which is canonically associated with
this class of ODEs. The conformal metric $g_{\cal N}$ yields all the
point 
invariant
information about the corresponding class of 3rd order ODEs.
\et

\section{Second order ODEs considered modulo point transformations}
This case has been recently carefully studied in
Ref. \cite{nurspar}. The ODE part of this paper includes, in
particular, description of the geometry associated with an equation 
\be
y''=Q(x,y,y')\label{ode2}
\ee
considered modulo point transformations (\ref{point}). This geometry,
in the convenient parametrization $(x,y,p=y')$ of the first jet space
$J^1$, turns out to be very closely related to the geometry associated
with the following split signature metric, the Fefferman metric,   
\be
g=2~[~(\der p-Q\der x)~\der x-(\der y-p\der x)~(\der\phi+\frac{2}{3}Q_p\der
  x+\frac{1}{6}Q_{pp}(\der y-p\der x))~]\label{feffer}
\ee
on $J^1\times{\bf R}$. More precisely, we have
the following theorem.
\bt~\\

\noindent
1) Every second order ODE (\ref{ode2})
endows its corresponding space $J^1\times{\bf R}$ with an orientation
and with the Fefferman metric (\ref{feffer}).\\

\noindent
2) If the ODE undergoes a point transformation (\ref{point}) then
  its Fefferman metric transforms conformally.\\

\noindent
3) All the point invariants of a point equivalence class of ODEs
  (\ref{ode2}) are expressible in terms of the conformal invariants of
the associated conformal class of Fefferman metrics.\\

\noindent
4) The Fefferman metrics (\ref{feffer}) are very special among all the split
signature metrics on 4-manifolds. Their Weyl tensor has algebraic type 
$(N,N)$ in the Cartan-Petrov-Penrose classification
{\rm \cite{c,pen,Pen1,pet}}. Both, the selfdual $C^+$ and the antiselfdual $C^-$,
parts of it are expressible in terms of only {\rm one} component. 
In fact, $C^+$ is proportional to 
$$w_1=D^2Q_{pp}-4DQ_{py}-DQ_{pp}Q_p+4Q_pQ_{py}-3Q_{pp}Q_y+6Q_{yy}$$
and $C^-$ is proportional to
$$w_2=Q_{pppp},$$
where
$$D=\partial_x+p\partial_y+Q\partial_p.$$
Each of the conditions $w_1=0$ and $w_2=0$ is invariant under 
point transformations (\ref{point}).\\

\noindent
5)  Cartan normal conformal connection 
associated with any conformal class $[g]$ of Fefferman 
metrics is {\rm reducible} to a certain 
${\bf SL}(2+1,{\bf R})$ connection naturally defined on an
8-dimensional bundle over $J^1$ which, via Cartan's equivalence
method, is uniquely associated with the point equivalence class of
corresponding ODEs (\ref{ode2}). The curvature of this connection has
very simple form
$$
\Om\sim
\begin{pmatrix}
0&w_2&*\\
&&\\
0&0&w_1\\
&&\\
0&0&0
\end{pmatrix}.
$$
If $w_1=0$ or $w_2=0$ this connection can be further understood as a
Cartan {\rm normal projective} connection over a certain two dimensional
space $S$ equipped with a projective structure {\rm \cite{Newman}}. 
$S$ can be identified either with the solution space of the ODE (\ref{ode2}) in
the $w_1=0$ case, or with the solution space of its {\rm
  dual}\footnote{See e.g. \cite{nurspar} for the concept of dual
  second order ODEs.} ODE in the $w_2=0$ case.  
\et
\section{Equations $z'=F(x,y,y',y'',z)$, noncompact form of the
  exceptional group $G_2$ and conformal metrics of signature $(3,2)$}
\subsection{Equations with integral-free solutions}
Consider a differential equation of the form
\be
G(x,y,y',...,y^{(m)},z,z',...,z^{(k)})=0\label{hil1}
\ee
for real functions $y=y(x)$ and $z=z(x)$ of one real variable. In
this equation $G:{\bf R}^{m+k+3}\to {\bf R}$ and $y^{(r)}$, $z^{(q)}$
denote the
$r$th and the $q$th derivative of $y$ and $z$ with respect to $x$. In 1912 Hilbert 
\cite{hilbert} considered a subclass of equations (\ref{hil1}) which
he called {\it equations with integral-free solutions} 
(Germ. integrallose Aufloesungen). These equations are defined as
follows.
\bd
Equation (\ref{hil1}) has {\rm integral-free solutions} iff its general
solution can be written as
\beq
&x=x(t,w(t),w'(t),...,w^{(p)}(t))\nonumber\\
&y=y(t,w(t),w'(t),...,w^{(p)}(t))\nonumber\\
&z=z(t,w(t),w'(t),...,w^{(p)}(t)),\nonumber
\eeq
where $w=w(t)$ is an {\rm arbitrary} sufficiently smooth real function of
one real variable.
\ed
\noindent
As an example consider equation 
\be
z'=y.\label{triv}
\ee 
Clearly $x=t$, $z=w(t)$,
$y=w'(t)$ is its general solution, which shows that (\ref{triv}) is in
the Hilbert class of equations with integral-free solutions. Very simple
equation (\ref{triv}) belongs to the class of {\it first order Monge equations}
\be
z'=F(x,y,y',z),\label{mong}
\ee 
which are equations (\ref{hil1}) with unknowns of at most of the first order.\\

\noindent
Associated with each first order Monge equation (\ref{mong}) there is a
4-dimensional space $J$ parametrized by $(x,y,p,z)$ and two 1-forms
\beq
&\om^1=\der z-F(x,y,p,z)\der x\nonumber\\
&\om^2=\der y-p\der x.\nonumber
\eeq
Every solution of the Monge equation (\ref{mong}) is a curve 
$
c(t)=(x(t),y(t),p(t),z(t))
$
in $J$ on which the forms $\om^1$ and $\om^2$ vanish.\\

\noindent
Suppose now, that given a Monge equation (\ref{mong}), there exists a 
transformation of the associated variables $(x,y,p,z)$
\be
\bma
x\\
y\\
p\\
z
\ema\stackrel{\phi}{\to}\bma
\bar{x}\\
\bar{y}\\
\bar{p}\\
\bar{z}
\ema=\bma
\bar{x}(x,y,p,z)\\
\bar{y}(x,y,p,z)\\
\bar{p}(x,y,p,z)\\
\bar{z}(x,y,p,z)
\ema\label{transs}
\ee
such that 
\beq
&\der\bar{y}-\bar{p}\der\bar{x}~=~\al\om^1+\bet\om^2\label{ffor}\\
&\der\bar{p}-\bar{z}\der\bar{x}~=~\gamma\om^1+\delta\om^2,\nonumber
\eeq
with $\al$,$\bet$, $\ga$, $\delta$  functions on $J$ satisfying
$\Delta=\al\delta-\bet\gamma\neq 0$. In such case 
\beq
&\om^1=\Delta^{-1}~[~\delta (\der\bar{y}-\bar{p}\der\bar{x})-\bet(\der\bar{p}-\bar{z}\der\bar{x})~]\nonumber \\
&\om^2=\Delta^{-1}~[-\gamma(\der\bar{y}-\bar{p}\der\bar{x})+\al(\der\bar{p}-\bar{z}\der\bar{x})~].\nonumber
\eeq
Thus, taking 
\be
\bar{x}=t,~~~\bar{y}=w(t),~~~\bar{p}=w'(t),~~~\bar{z}=w''(t)\label{subst}
\ee
we construct a curve in $J$ on which the forms $\om^1$ and $\om^2$ 
identically vanish. Now, the inverse of $\phi$ which 
gives $x=x(\bar{x},\bar{y},\bar{p},\bar{z})$, etc., provides
\beq
&x=x(t,w(t),w'(t),w''(t))\nonumber\\
&y=y(t,w(t),w'(t),w''(t))\nonumber\\
&z=z(t,w(t),w'(t),w''(t)),\nonumber
\eeq
which is an integral-free solution of the Monge equation (\ref{mong}).\\

\noindent
We summarize our discussion in the following Lemma.
\bl~\\
Every first order Monge equation (\ref{mong}) admitting coordinate 
transformation (\ref{transs}) which realizes 
(\ref{ffor}) has integral-free solutions.
\el

\noindent
{\bf Example 4}\\
Consider equation
\be
z'=(y')^2.\label{exeq}
\ee
Its corresponding forms are $\om^1=\der z-p^2\der x$, $\om^2=\der
y-p\der x$. The change of variables $x=\frac{1}{2}\bar{z}$,
$y=\frac{1}{2}(\bar{z}\bar{x}-\bar{p})$,~ 
$z=\frac{1}{2}\bar{z}\bar{x}^2-\bar{p}\bar{x}+\bar{y}$, ~$p=\bar{x}$ 
brings them to the form
$\om^1=\der\bar{y}-\bar{p}\der\bar{x}-\bar{x}(\der\bar{p}-
\bar{z}\der\bar{x})$,
~$\om^2=-\frac{1}{2}(\der\bar{p}-\bar{z}\der\bar{x})$. 
This proves that substitution (\ref{subst}) leads to the following integral-free solution of equation (\ref{exeq}):
\beq
&x=\frac{1}{2}w''(t)\nonumber\\
&y=\frac{1}{2}t w''(t)-\frac{1}{2}w'(t)\nonumber\\
&z=\frac{1}{2}t^2w''(t)-tw'(t)+w(t).\nonumber
\eeq 

\noindent
A natural question as to whether all the first order 
Monge equations have integral-free
 solutions was answered in affirmative by Monge. Thus, we
 have the following theorem.

\bt (Monge)\\
Every first order Monge equation has integral-free solutions. 
\et

\noindent
It is instructive to sketch the proof of this theorem.\\

\noindent
Given a Monge equation (\ref{mong}) we consider its associated two
1-forms
\be
\om^1=\der z-F(x,y,p,z)\der x~~~~~~~{\rm and}~~~~~~~\om^2=\der y
-p\der x\label{fff}
\ee
on $J$. We say that another pair 
of linearly independent 1-forms $(\bar{\om}^1,\bar{\om}^2)$ on $J$ 
is {\it equivalent} to the pair (\ref{fff}) if there exists a
transformation of variables (\ref{transs}) and functions $\al,\bet,\ga,\delta$,
$\al\delta-\bet\ga\neq 0$, on $J$ such that  
\beq
&\phi^*(\bar{\om}^1)~=~\al\om^1+\bet\om^2\label{ffop}\\
&\phi^*(\bar{\om}^2)~=~\gamma\om^1+\delta\om^2.\nonumber
\eeq 
\noindent
According to Lemma 1, if we were able to show that there is only one 
equivalence class of forms $(\om^1,\om^2)$ equivalent to 
$(\der \bar{y}-\bar{p}\der
\bar{x},\der\bar{p}-\bar{z}\der\bar{x})$, the theorem would be
proven. Thus, in the process of proving the Monge theorem, we are led to
study the equivalence problem for two 1-forms given modulo
transformations (\ref{ffop})  on an open set of ${\bf
  R}^4$. Introducing the total differential
vector field $D=\partial_x+p\partial_y+F\partial_z$ it is not difficult
to prove that
a pair of 1-forms (\ref{fff}) originating from the 
Monge equations for which 
\be
F_{pp}=0 ~~~~~{\rm and}~~~~~ DF_p-F_y-F_pF_z=0\label{ccon}
\ee
and a pair of forms originating from the equations for which at least one of the above conditions is not
satisfied are {\it not} equivalent. Then, the Cartan equivalence
method applied to the forms related to the first order 
Monge equations {\it not}
satisfying (\ref{ccon}) shows that they are {\it all} locally equivalent to 
$(\der \bar{y}-\bar{p}\der
\bar{x},\der\bar{p}-\bar{z}\der\bar{x})$. Thus, the first order 
Monge equations 
for which at least one of conditions (\ref{ccon}) is not satisfied
have general solutions of the form 
\beq
&x=x(t,w(t),w'(t),w''(t))\nonumber\\
&y=y(t,w(t),w'(t),w''(t))\label{cc1}\\
&z=z(t,w(t),w'(t),w''(t)).\nonumber
\eeq
On the other hand, if we apply the Cartan equivalence method to the
forms originating from the Monge equations {\it satisfying} (\ref{ccon}), we
show that they are {\it all} locally equivalent
to  $(\der \bar{z},\der\bar{y}-\bar{p}\der\bar{x})$. Thus, taking
$\bar{z}=$const, $\bar{x}=t$, $\bar{y}=w(t)$ and $\bar{p}=w'(t)$ we
show that in such case the Monge equations have general solutions of
the form 
\beq
&x=x(t,w(t),w'(t))\nonumber\\
&y=y(t,w(t),w'(t))\label{cc2}\\
&z=z(t,w(t),w'(t)).\nonumber
\eeq
Therefore in the both nonequivalent cases (\ref{cc1}) and (\ref{cc2}) the 
Monge equations have integral-free solutions. This finishes the proof
of the Monge theorem.\\

\noindent
Hilbert in \cite{hilbert} considered an equation 
\be
z'=(y'')^2\label{hill}
\ee
and proved that it has {\it not} the property of having
integral-free solutions. It turns out, that among all the equations
which have not this property, the Hilbert equation (\ref{hill}) is, in
a certain sense, the simplest one. 
\subsection{Equivalence of forms associated with ODEs  $z'=F(x,y,y',y'',z)$}
The Hilbert equation (\ref{hill}) is a special case of an equation 
\be
z'=F(x,y,y',y'',z).\label{hilb}
\ee
Equations of this type were considered by Cartan \cite{5var} who, in
particular, 
observed that they describe 
Cauchy characteristics of pairs of involutive second order PDEs for a real
function of two variables. In the context of the present paper we
are interested under what conditions equations (\ref{hilb}) have
integral-free solutions. The treatment of the problem is a simple
generalization of the method described in the sketch of the proof of
Monge's theorem. Thus, with each equation (\ref{hilb}) we associate
{\it three} 1-forms 
\beq
&\om^1=\der z-F(x,y,p,q,z)\der x\nonumber\\
&\om^2=\der y-p\der x\label{fhil}\\
&\om^3=\der p-q\der x,\nonumber
\eeq
which live on a 5-dimensional manifold $J$ parametrized by
$(x,y,p=y',q=y'',z)$. Following the case of Monge equations, we need
to study the equivalence problem for the {\it triples} of linearly
independent 1-forms $(\om^1,\om^2,\om^3)$ on an open set of ${\bf
  R}^5$. More precisely, let $(\om^1,\om^2,\om^3)$ be defined on a
open set $J\subset {\bf R}^5$ 
parametrized by $(x,y,p,q,z)$ and
$(\bar{\om}^1,\bar{\om}^2,\bar{\om}^3)$ be defined on a set
$\bar{J}\subset {\bf R}^5$
parametrized by $(\bar{x},\bar{y},\bar{p},\bar{q},\bar{z})$. We say
that the two triples $(\om^1,\om^2,\om^3)$ and
$(\bar{\om}^1,\bar{\om}^2,\bar{\om}^3)$ are (locally) 
equivalent iff there exists a (local) diffeomorphism $\phi:J\to\bar{J}$ 
\be
\bma
x\\
y\\
p\\
q\\
z
\ema\stackrel{\phi}{\to}\bma
\bar{x}\\
\bar{y}\\
\bar{p}\\
\bar{q}\\
\bar{z}
\ema=\bma
\bar{x}(x,y,p,q,z)\\
\bar{y}(x,y,p,q,z)\\
\bar{p}(x,y,p,q,z)\\
\bar{q}(x,y,p,q,z)\\
\bar{z}(x,y,p,q,z)
\ema\label{transss}
\ee
and a
${\bf GL}(3,{\bf R})$-valued function 
$$
f=\bma
\al&\bet&\ga\\
\delta&\epsilon&\lambda\\
\kappa&\mu&\nu\\
\ema
$$
on $J$ such that 
\beq
&\phi^*(\bar{\om}^1)~=~\al\om^1+\bet\om^2+\ga\om^3\nonumber\\
&\phi^*(\bar{\om}^2)~=~\delta\om^1+\epsilon\om^2+\lambda\om^3\label{transp}\\
&\phi^*(\bar{\om}^3)~=~\kappa\om^1+\mu\om^2+\nu\om^3.\nonumber
\eeq
The equivalence problem for such triples was solved
by Cartan. His solution, in particular, can be applied to the triples
of 1-forms (\ref{fhil}) originating from the Cartan equations
(\ref{hilb}). Cartan's analysis, restricted to such triples, shows
that they split onto two main nonequivalent classes. The first class 
originates from equations (\ref{hilb}) satisfying 
$$
F_{qq}=0,
$$
the second class is defined by the equations for which 
$$
F_{qq}\neq 0.
$$ 
Both the above classes include nonequivalent triples of 1-forms, but
only the first class originates from equations (\ref{hilb}) with integral-free
solution. All the Cartan equations with $F_{qq}\neq 0$ have not the
property of having integral-free solutions. The Hilbert equation
(\ref{hill}) is one of the equations from this class.\\

\noindent
{\bf Example 5}\\
According to the above discussion, if $k\neq 0$ and $k\neq 1$ equation 
\be
z'=\frac{1}{k}(y'')^k\label{moj}
\ee 
has not the property of having integral-free
solutions. Thus, since one is forced to use integrals to write down
the general solution of (\ref{moj}), we solve it by putting 
\be
x=t,~~~~~y=w(t),~~~~~~z=\frac{1}{k}\int w''(t)^k\der t.\label{brut}
\ee
Cartan found {\it better} solution
\beq
&x=(k-1)t^{\frac{k-2}{k-1}}w''(t)\nonumber\\
&\nonumber\\
&y=\frac{1}{2}(k-1)^2t^{\frac{2k-3}{k-1}}w''(t)^2-(k-1)t^{\frac{k-2}{k-1}}w'(t)w''(t)+\frac{1}{2}(k-1)\int
t^{\frac{k-2}{k-1}}w''(t)^2\der t\nonumber\\
&\nonumber\\
&z=\frac{k-1}{k}t^2w''(t)-tw'(t)+w(t).\nonumber
\eeq   
We prefer this solution rather then (\ref{brut}) since it involves
only second power of $w''$ under the integral, whereas the solution
(\ref{brut}) involves the $k$th power. This example shows that, for a
given Cartan equation, among
many different expressions for its general solution which involve
integrals there could be some preferred ones. The precise meaning of
this observation is worth further investigation.  \\

\subsection{$\tilde{G}_2$ Cartan connection for equation
  $z'=F(x,y,y',y'',z)$ and conformal (3,2)-signature geometry } 
We will not comment any further on Cartan equations
for which $F_{qq}=0$. Instead, we concentrate on much more
interesting $F_{qq}\neq 0$ case. \\

\noindent
First, we briefly sketch Cartan's results on equivalence problem for
forms 
\beq
&\om^1=\der z-F(x,y,p,q,z)\der x\nonumber\\
&\om^2=\der y-p\der x\label{55p}\\
&\om^3=\der p-q\der x\nonumber
\eeq
satisfying $F_{qq}\neq 0$. On doing that we supplement these forms to 
a coframe $(\om^1,\om^2,\om^3,\om^4,\om^5)$ on the $(x,y,p,q,z)$ space
such that 
\beq
&\om^1=\der z-F(x,y,p,q,z)\der x\nonumber\\
&\om^2=\der y-p\der x\nonumber\\
&\om^3=\der p-q\der x\label{com}\\
&\om^4=\der q\nonumber\\
&\om^5=\der x.\nonumber
\eeq
Since we are interested in {\it all} forms $(\om^1,\om^2,\om^3)$ which
are equivalent
to the forms (\ref{55p}) via transformations
(\ref{transss})-(\ref{transp}) this coframe is not
unique. It is given up to the following freedom:
$$
\bma
\om^1\\
\om^2\\
\om^3\\
\om^4\\
\om^5
\ema\to\bma
\bar{\om}^1\\
\bar{\om}^2\\
\bar{\om}^3\\
\bar{\om}^4\\
\bar{\om^5}
\ema=\bma
\al&\bet&\ga&0&0\\
\delta&\epsilon&\lambda&0&0\\
\kappa&\mu&\nu&0&0\\
\pi&\rho&\si&\tau&\chi\\
\pi'&\rho'&\si'&\tau'&\chi'
\ema
\bma
\om^1\\
\om^2\\
\om^3\\
\om^4\\
\om^5
\ema,
$$   
which suggests that instead of working with a not uniquely defined
coframe (\ref{com}) on the $(x,y,p,q,z)$ space it is better to use
five well defined linearly independent 1-forms 
$$
\bma
\theta^1\\
\theta^2\\
\theta^3\\
\theta^4\\
\theta^5
\ema =\bma
\al&\bet&\ga&0&0\\
\delta&\epsilon&\lambda&0&0\\
\kappa&\mu&\nu&0&0\\
\pi&\rho&\si&\tau&\chi\\
\pi'&\rho'&\si'&\tau'&\chi'
\ema
\bma
\om^1\\
\om^2\\
\om^3\\
\om^4\\
\om^5
\ema,
$$  
on a bigger space parametrized by
$(x,y,p,q,z,\al,\bet,\ga,\delta,\epsilon,\lambda,\kappa,\mu,\nu,\pi,\rho,\sigma,\tau,\chi,\pi',\rho',\sigma',\tau',\chi')$.
Now, assuming that $F_{qq}\neq 0$ and using his equivalence method (which
involved several reductions and prolongations\footnote{See e.g. in
  Ref. \cite{Olver} for the definitions of these procedures}) Cartan
was able to prove that on a certain 14-dimensional manifold $P$ the
forms $(\theta^1,\theta^2,\theta^3,\theta^4,\theta^5)$ can be
supplemented in a unique way to a unique coframe. More precisely, he
proved the following theorem.
\bt~(Cartan)\\
An equivalence class of forms 
\beq
&\om^1=\der z-F(x,y,p,q,z)\der x\nonumber\\
&\om^2=\der y-p\der x\label{cof}\\
&\om^3=\der p-q\der x,\nonumber
\eeq
for which $F_{qq}\neq 0$, {\rm uniquely} defines a
14-dimensional manifold $P$ and a preferred coframe
$(\theta^1,\theta^2,\theta^3,$ 
$\theta^4,\theta^5,\Om_1,\Om_2,\Om_3,\Om_4,\Om_5,\Om_6,\Om_7,\Om_8,\Om_9)$
on it such that
\beq
&\der\theta^1=\theta^1\dz(2\Om_1+\Om_4)+\theta^2\dz\Om_2+\theta^3\dz\theta^4\nonumber\\
&\nonumber\\
&\der\theta^2=\theta^1\dz\Om_3+\theta^2\dz(\Om_1+2\Om_4)+\theta^3\dz\theta^5\nonumber\\&\nonumber\\
&\der\theta^3=\theta^1\dz\Om_5+\theta^2\dz\Om_6+\theta^3\dz(\Om_1+\Om_4)+\theta^4\dz\theta^5\label{sycart}\\&\nonumber\\
&\der\theta^4=\theta^1\dz\Om_7+\frac{4}{3}\theta^3\dz\Om_6+\theta^4\dz\Om_1+\theta^5\dz\Om_2\nonumber\\&\nonumber\\
&\der\theta^5=\theta^2\dz\Om_7-\frac{4}{3}\theta^3\dz\Om_5+\theta^4\dz\Om_3+\theta^5\dz\Om_4.\nonumber
\eeq
\et  
Note that the above theorem implies formulae for the differentials of
the forms $\Om_\mu$,
$\mu=1,2,...,9$. Explicitly, these differentials are:
\beq
&\der\Om_1=\Om_3\dz\Om_2+\frac{1}{3}\theta^3\dz\Om_7-\frac{2}{3}\theta^4\dz\Om_5+\frac{1}{3}\theta^5\dz\Om_6+\theta^1\dz\Om_8+\nonumber\\
&\frac{3}{8}c_2\theta^1\dz\theta^2+b_2\theta^1\dz\theta^3+b_3\theta^2\dz\theta^3+a_2\theta^1\dz\theta^4+a_3\theta^1\dz\theta^5+a_3\theta^2\dz\theta^4+a_4\theta^2\dz\theta^5\nonumber\\
&\nonumber\\
&\der\Om_2=\Om_2\dz(\Om_1-\Om_4)-\theta^4\dz\Om_6+\theta^1\dz\Om_9+\nonumber\\
&\frac{3}{8}c_3\theta^1\dz\theta^2+b_3\theta^1\dz\theta^3+a_3\theta^1\dz\theta^4+a_4\theta^1\dz\theta^5+b_4\theta^2\dz\theta^3+a_4\theta^2\dz\theta^4+a_5\theta^2\dz\theta^5\nonumber\\&\nonumber\\
&\der\Om_3=\Om_3\dz(\Om_4-\Om_1)-\theta^5\dz\Om_5+\theta^2\dz\Om_8-\nonumber\\
&\frac{3}{8}c_1\theta^1\dz\theta^2-b_1\theta^1\dz\theta^3-a_1\theta^1\dz\theta^4-a_2\theta^1\dz\theta^5-b_2\theta^2\dz\theta^3-a_2\theta^2\dz\theta^4-a_3\theta^2\dz\theta^5\nonumber\\&\nonumber\\
&\der\Om_4=\Om_2\dz\Om_3+\frac{1}{3}\theta^3\dz\Om_7+\frac{1}{3}\theta^4\dz\Om_5-\frac{2}{3}\theta^5\dz\Om_6+\theta^2\dz\Om_9-\nonumber\\
&\frac{3}{8}c_2\theta^1\dz\theta^2-b_2\theta^1\dz\theta^3-a_2\theta^1\dz\theta^4-a_3\theta^1\dz\theta^5-b_3\theta^2\dz\theta^3-a_3\theta^2\dz\theta^4-a_4\theta^2\dz\theta^5\nonumber\\&\nonumber\\
&\der\Om_5=\Om_1\dz\Om_5+\Om_3\dz\Om_6-\theta^5\dz\Om_7+\theta^3\dz\Om_8+\label{syspp}\\
&\frac{9}{32}\delta_1\theta^1\dz\theta^2+\frac{3}{4}c_1\theta^1\dz\theta^3+\frac{3}{4}b_1\theta^1\dz\theta^4+\frac{3}{4}b_2\theta^1\dz\theta^5+\frac{3}{4}c_2\theta^2\dz\theta^3+\frac{3}{4}b_2\theta^2\dz\theta^4+\frac{3}{4}b_3\theta^2\dz\theta^5\nonumber\\
&\nonumber\\
&\der\Om_6=\Om_2\dz\Om_5+\Om_4\dz\Om_6+\theta^4\dz\Om_7+\theta^3\dz\Om_9+\nonumber\\
&\frac{9}{32}\delta_2\theta^1\dz\theta^2+\frac{3}{4}c_2\theta^1\dz\theta^3+\frac{3}{4}b_2\theta^1\dz\theta^4+\frac{3}{4}b_3\theta^1\dz\theta^5+\frac{3}{4}
c_3\theta^2\dz\theta^3+\frac{3}{4}b_3\theta^2\dz\theta^4+\frac{3}{4}b_4\theta^2\dz\theta^5\nonumber\\&\nonumber\\
&\der\Om_7=\frac{4}{3}\Om_5\dz\Om_6+(\Om_1+\Om_4)\dz\Om_7+\theta^4\dz\Om_8+\theta^5\dz\Om_9+\nonumber\\
&\frac{9}{64}e\theta^1\dz\theta^2-\frac{3}{8}\delta_1\theta^1\dz\theta^3-\frac{3}{8}c_1\theta^1\dz\theta^4-\frac{3}{8}c_2\theta^1\dz\theta^5-\frac{3}{8}\delta_2\theta^2\dz\theta^3-\frac{3}{8}c_2\theta^2\dz\theta^4-\frac{3}{8}c_3\theta^2\dz\theta^5\nonumber\\&\nonumber\\
&\der\Om_8=\Om_5\dz\Om_7+(2\Om_1+\Om_4)\dz\Om_8+\Om_3\dz\Om_9+\nonumber\\
& h_1\theta^1\dz\theta^2 +
h_2\theta^1\dz\theta^3 + h_3\theta^1\dz\theta^4 +
h_4\theta^1\dz\theta^5 + h_5\theta^2\dz\theta^3 + h_4\theta^2\dz\theta^4 +
h_6\theta^2\dz\theta^5\nonumber\\
&\nonumber\\
&\der\Om_9=\Om_6\dz\Om_7+(\Om_1+2\Om_4)\dz\Om_9+\Om_2\dz\Om_8+\nonumber\\
& k_1\theta^1\dz\theta^2 +\frac{1}{32}(3e+32h_5)\theta^1\dz\theta^3 + \frac{1}{32}(-3 \delta_1
  + 32 h_4)\theta^1\dz\theta^4 +\nonumber\\
&\frac{1}{32}(-3\delta_2+32h_6)\theta^1\dz\theta^5 + k_2\theta^2\dz\theta^3 + \frac{1}{32}(-3\delta_2 +32h_6)\theta^2\dz\theta^4 + k_3\theta^2\dz\theta^5\nonumber
\eeq 
where $a_1$, $a_2$, $a_3$, $a_4$, $a_5$, $b_1$, $b_2$, $b_3$, $b_4$,
$c_1$, $c_2$, $c_3$, $\delta_1$, $\delta_2$, $e$, $h_1$, $h_2$, $h_3$, $h_4$, $h_5$, $h_6$, $k_1$, $k_2$, $k_3$ are functions on $P$
uniquely defined by the equivalence class of forms (\ref{cof}).\\

\noindent
The system (\ref{sycart})-(\ref{syspp}) provides all the local 
invariants for the equivalence class of forms (\ref{cof}) satisfying
$F_{qq}\neq 0$. If one is given two triples of
1-forms 
$$\om^1=\der z-F(x,y,p,q,z)\der x,~~~~~~~~F_{qq}\neq 0,$$
$$\om^2=\der y-p\der x$$
$$\om^3=\der p-q\der x$$
and 
$$\bar{\om}^1=\der
\bar{z}-\bar{F}(\bar{x},\bar{y},\bar{p},\bar{q},\bar{z})
\der \bar{x},~~~~~~~~\bar{F}_{\bar{q}\bar{q}}\neq 0,$$
$$\bar{\om}^2=\der \bar{y}-\bar{p}\der \bar{x}$$
$$\bar{\om}^3=\der \bar{p}-\bar{q}\der \bar{x}$$
on respective manifolds $J$ and $\bar{J}$ parametrized by
$(x,y,p,q,z)$ and $(\bar{x},\bar{y},\bar{p},\bar{q},\bar{z})$, then
there exists a local diffeomorphism 
$$
\bma
x\\
y\\
p\\
q\\
z
\ema\stackrel{\phi}{\to}\bma
\bar{x}\\
\bar{y}\\
\bar{p}\\
\bar{q}\\
\bar{z}
\ema=\bma
\bar{x}(x,y,p,q,z)\\
\bar{y}(x,y,p,q,z)\\
\bar{p}(x,y,p,q,z)\\
\bar{q}(x,y,p,q,z)\\
\bar{z}(x,y,p,q,z)\ema
$$
realizing 
$$\phi^*(\bar{\om}^1)~=~\al\om^1+\bet\om^2+\ga\om^3$$
$$\phi^*(\bar{\om}^2)~=~\delta\om^1+\epsilon\om^2+\lambda\om^3$$
$$\phi^*(\bar{\om}^3)~=~\kappa\om^1+\mu\om^2+\nu\om^3$$
iff there exists a diffeomorphism $\Phi:P\to\bar{P}$ between the
associated 14-dimensional manifolds $P$ and $\bar{P}$ of Theorem 8 such that 
$$\Phi^*(\bar{\theta}^i)=\theta^i,~~~~~~~~~\Phi^*(\bar{\Om}_\mu)=\Om_\mu$$
for all $i=1,2,3,4,5$ and $\mu=1,2,3,...,9$. This, in particular means
that to realize the equivalence between the $(\om^i)$s and
$(\bar{\om}^i)$s, the diffeomorphism $\Phi$ must also satisfy 
$$\Phi^*(\bar{a}_1)=a_1,~~~\Phi^*(\bar{b}_1)=b_1,~~~\Phi^*(\bar{c}_1)=c_1,~~~
{\rm etc.}$$
This gives severe algebraic (i.e. non-differential) constraints on
$\Phi$ and, in generic cases, quickly leads to the answer if the two
systems of forms $(\om^i)$ and $(\bar{\om}^i)$ are equivalent.\\

\noindent
In view of the above we ask for those equivalence classes of forms
(\ref{cof}) which correspond to systems (\ref{sycart})-(\ref{syspp})
with all the {\it scalar invariants}
$(a_1,a_2,a_3,a_4,a_5,b_1,b_2,b_3,b_4,c_1,c_2,c_3,\delta_1,\delta_2,e,h_1,h_2,h_3,h_4,$
$h_5,h_6,k_1,k_2,k_3)$ being
constants. It follows that it is possible if and only if all of them
are identically equal to zero. In this well defined case the system 
(\ref{sycart})-(\ref{syspp}) can be understood as a system consisting of 
right invariants forms $(\theta^i,\Om_\mu)$ on a 14-dimensional 
Lie group. This group is simple and has indefinite Killing form, as
can be seen from the structure constant coefficients defined by the
system (\ref{sycart})-(\ref{syspp}) with all the scalar invariants
vanishing. This identifies this group as a
noncompact real form $\tilde{G}_2$ of the exceptional group
$G_2$. \\

\noindent
It follows that there is only one equivalence class of forms (\ref{cof}) 
corresponding to the system (\ref{sycart})-(\ref{syspp}) with all the 
scalar invariants vanishing. It can be defined by the function 
$$F=q^2$$ associated with the Hilbert equation
$$z'=(y'')^2.$$\\

\noindent
In case of general scalar invariants, the system
(\ref{sycart})-(\ref{syspp}) defines a curvature of a certain Cartan
$\underline{\tilde{\bf g}_2}$-valued connection which `measures' how much
the equivalence class of forms (\ref{cof}) is distorted from the flat
Hilbert case corresponding to $F=q^2$. To define this connection we
first observe that the system (\ref{sycart})-(\ref{syspp}) guarantees
that $P$ is foliated by 9-dimensional leaves. These are the integral
manifolds of the distribution spanned by 
vector fields $Y_\mu$, $\mu=1,2,...9$ which, together with $X_i$,
$i=1,2,...5$, form a frame 
$(X_1,X_2,X_3,X_4,X_5,Y_1,Y_2,Y_3,Y_4,Y_5,Y_6,Y_7,Y_8,Y_9)$ 
dual to the invariant coframe 
$(\theta^1,\theta^2,\theta^3,\theta^4,\theta^5,\Om_1,\Om_2,\Om_3,\Om_4,\Om_5,\Om_6,\Om_7,\Om_8,\Om_9)$
on $P$. (The fact that this distribution is integrable, is a
simple corollary, from equations (\ref{sycart}), which show that the
basis $\theta^i$, $i=1,2,..5$, of its annihilator is a differential
ideal). This proves that the manifold $P$ is fibered over a 5-dimensional
space of leaves of this distribution. This space may be identified
with the $(x,y,p,q,z)$ space $J$ on which the original forms $\om^i$,
$i=1,2,...5$, defining the equivalence class (\ref{cof}) reside. Thus
we have a fibration $P\to J$, which is actually a principal fibre
bundle with the 9-dimensional parabolic subgroup $H$ of $\tilde{G}_2$
as its structure group. On this fibre bundle the following matrix of 1-forms:
\be
\om_{\tilde{G}_2}=\bma
-\Om_1-\Om_4&-\Om_8&-\Om_9&-\frac{1}{\sqrt{3}}\Om_7&\frac{1}{3}\Om_5&\frac{1}{3}\Om_6&0\\&&&&&&\\
\theta^1&\Om_1&\Om_2&\frac{1}{\sqrt{3}}\theta^4&-\frac{1}{3}\theta^3&0&\frac{1}{3}\Om_6\\&&&&&&\\
\theta^2&\Om_3&\Om_4&\frac{1}{\sqrt{3}}\theta^5&0& -\frac{1}{3}\theta^3&-\frac{1}{3}\Om_5\\&&&&&&\\
\frac{2}{\sqrt{3}}\theta^3&\frac{2}{\sqrt{3}}\Om_5& \frac{2}{\sqrt{3}}\Om_6&0&
\frac{1}{\sqrt{3}}\theta^5&-\frac{1}{\sqrt{3}}\theta^4&-\frac{1}{\sqrt{3}}\Om_7\\&&&&&&\\
\theta^4& \Om_7&0&\frac{2}{\sqrt{3}}\Om_6&-\Om_4& \Om_2&
\Om_9\\&&&&&&\\
\theta^5& 0& \Om_7&-\frac{2}{\sqrt{3}}\Om_5& \Om_3& -\Om_1&-\Om_8\\&&&&&&&\\
0&\theta^5&-\theta^4&
\frac{2}{\sqrt{3}}\theta^3&-\theta^2& \theta^1&\Om_1+\Om_4
\ema,\label{ccg2}
\ee  
becomes a Cartan connection with values in the Lie algebra of
$\tilde{G}_2$. (The fact that $\om_{\tilde{G}_2}$ is 
$\underline{\tilde{\bf g}_2}$-valued can be checked e.g. by successive replacement of one of
the 14 forms
$(\theta^1,\theta^2,\theta^3,\theta^4,\theta^5,\Om_1,$ $\Om_2,\Om_3,\Om_4,\Om_5,\Om_6,\Om_7,\Om_8,\Om_9)$
in $\om_{\tilde{G}_2}$ by 1 with simultaneous replacement of all the
  others forms by 0.  The so obtained 14 matrices satisfy the commutation
  relations of $\underline{\tilde{\bf g}_2}$.) The curvature of
  this connection 
$${\cal R}=\der\om_{\tilde{G}_2}+\om_{\tilde{G}_2}\dz\om_{\tilde{G}_2},$$
being horizontal, involves only $\theta^i\dz\theta^j$ terms. This when
compared with equations (\ref{syspp}), enables the scalar invariants
to be interpreted as the curvature coefficients of $\om_{\tilde{G}_2}$.\\

\noindent
Another interpretation of $\om_{\tilde{G}_2}$ can be obtained by
recalling that $\tilde{G}_2$ is naturally embedded in ${\bf SO}(4,3)$
as its subgroup stabilizing a generic 3-form in ${\bf
  R}^{(4,3)}$. We have chosen a 7-dimensional representation
of the Lie algebra $\underline{\tilde{\bf g}_2}$ in such a way that
the connection 
$\om_{\tilde{G}_2}$ can be interpreted as a reduction of a Cartan
normal conformal connection associated with a certain
$(3,2)$-signature conformal metric defined on $J$. In the following we
describe this view point.\\

\noindent
Given an equivalence class of forms (\ref{cof}) satisfying $F_{qq}\neq 0$
and using the forms 
$(\theta^1,\theta^2,\theta^3,\theta^4,\theta^5)$ associated with them
via Theorem 8 we define a following bilinear form 
\be
\tilde{g}=2\theta^1\theta^5-2\theta^2\theta^4
+\frac{4}{3}\theta^3\theta^3\label{met32}
\ee
on $P$. This form is clearly degenerate and has signature $(+,+,+,-,-,0,0,0,0,0,0,0,0,0)$. Using the frame
$(X_1,X_2,X_3,X_4,X_5,Y_1,Y_2,Y_3,Y_4,Y_5,Y_6,Y_7,Y_8,Y_9)$ on $P$ defined
above, we see that the degenerate directions of $\tilde{g}$ are tangent to the
vectors $Y_\mu$. Now, the system (\ref{sycart}) guarantees that 
the form $\tilde{g}$ scales when Lie dragged along any of the 
directions $Y_\mu$. In other words we have 
$${\cal L}_{Y_\mu}\tilde{g}=\lambda_\mu\tilde{g}$$
with some functions $\lambda_\mu$. This, when compared with the fact
that the distribution spanned by $Y_\mu$, $\mu=1,2,...9$, 
defines a foliation on $P$, means that 
the degenerate bilinear form $\tilde{g}$ projects from $P$ to $J$, the
space of leaves of this foliation, defining there a conformal metric 
$[G_{(3,2)}]$ of signature $(+,+,+,-,-)$. It is this conformal
structure that yields all the information about the local invariants
of an equivalence class of forms (\ref{cof}). Calculating the Cartan
normal conformal connection of this conformal structure, leads to the 
conclusion that it is reducible to the $\underline{\tilde{\bf g}_2}$-valued 
Cartan connection $\om_{\tilde{G}_2}$ on $P$.\\

\noindent
Remarkably the conformal metric $[G_{(3,2)}]$ is defined on the same
space $J$ on which the original forms $\om^i$, $i=1,2,...5$, defining
the equivalence class (\ref{cof}) were defined. Thus, it is possible
to write down a local representative $G_{(3,2)}$ of $[G_{(3,2)}]$ in 
coordinates $(x,y,p,q,z)$ in which the forms $\om^i$ read
\beq
&\om^1=\der z-F(x,y,p,q,z)\der x\nonumber\\
&\om^2=\der y-p\der x\nonumber\\
&\om^3=\der p-q\der x\nonumber\\
&\om^4=\der q\nonumber\\
&\om^5=\der x.\nonumber
\eeq
Introducing the total differential operator $D$ on $J$ by 
$$
D=\partial_x+p\partial_y+q\partial_p+F\partial_z
$$     
we find that a representative of $[G_{(3,2)}]$ is given by
\beq
&G_{(3,2)}=\label{32met}\\
&[~DF_{qq}^2 F_{qq}^2 + 6DF_q DF_{qqq}F_{qq}^2 - 6DF_{qqq}F_p F_{qq}^2
  -3DDF_{qq} F_{qq}^3 + 9 DF_{qp}F_{qq}^3 - 9 F_{pp}F_{qq}^3 +\nonumber\\&9
  DF_{qz} F_q F_{qq}^3 - 18 F_{pz} F_qF_{qq}^3 + 3DF_z
  F_{qq}^4 - 6 DF_qF_{qq}^2 F_{qqp} + 6 F_p F_{qq}^2 F_{qqp} - \nonumber\\&8
  DF_qDF_{qq}F_{qq}F_{qqq} + 8 DF_{qq}F_p F_{qq}F_{qqq} +
  3DDF_qF_{qq}^2 F_{qqq} - 3DF_p F_{qq}^2 F_{qqq} - 3DF_z F_q F_{qq}^2
  F_{qqq} +\nonumber\\& 4(DF_q)^2 F_{qqq}^2 - 8 DF_q F_pF_{qqq}^2  - 3 (DF_q)^2
  F_{qq}F_{qqqq}+ 4 F_p^2 F_{qqq}^2+ 6DF_q F_p
  F_{qq}F_{qqqq} -\nonumber\\&  3 F_p^2 F_{qq}F_{qqqq} - 6 DF_q F_q F_{qq}^2
  F_{qqz} +6 F_p F_q F_{qq}^2 F_{qqz} - 3 DF_q F_{qq}^3 F_{qz}
  + 12 F_p F_{qq}^3 F_{qz} +\nonumber\\& 3 F_{qq}^2 F_{qqq}F_y - 6
  DF_{qqq}F_qF_{qq}^2 F_z + 4DF_{qq} F_{qq}^3 F_z  + 6 F_q F_{qq}^2
  F_{qqp}F_z + 8 DF_{qq}F_qF_{qq}F_{qqq}F_z -\nonumber\\& 4 DF_q
  F_{qq}^2 F_{qqq} F_z- 9 F_{qp} F_{qq}^3 F_z + F_p F_{qq}^2 F_{qqq}
  F_z - 8 DF_q F_qF_{qqq}^2 F_z + 8 F_p F_q F_{qqq}^2 F_z
  +\nonumber\\&
 6 DF_q F_q
  F_{qq}F_{qqqq} F_z - 6 F_p F_q F_{qq} F_{qqqq} F_z + 18
  F_{qq}^3 F_{qy} + 6 F_q^2 F_{qq}^2 F_{qqz} F_z + 3F_q F_{qq}^3
  F_{qz}F_z -\nonumber\\& 2 F_{qq}^4 F_z^2 + F_q F_{qq}^2 F_{qqq} F_z^2 + 4 F_q^2
  F_{qqq}^2 F_z^2- 3 F_q^2 F_{qq} F_{qqqq}F_z^2 - 9 F_q^2
  F_{qq}^3 F_{zz}~]~(\tilde{\om}^1)^2+\nonumber\\
&[~6 DF_{qqq} F_{qq}^2 - 6 F_{qq}^2 F_{qqp} - 8
  DF_{qq}F_{qq}F_{qqq} + 
8 DF_q F_{qqq}^2 - 8 F_p F_{qqq}^2 - 6 DF_q F_{qq} F_{qqqq} + 
\nonumber\\&6 F_p F_{qq}F_{qqqq} - 6 F_q F_{qq}^2F_{qqz} + 6 F_{qq}^3 F_{qz} + 2
F_{qq}^2 F_{qqq} F_z - 8 F_q F_{qqq}^2 F_z + 6 F_q F_{qq}F_{qqqq}F_z~]~\tilde{\om}^1\tilde{\om}^2 +\nonumber\\& [~10 DF_{qq}F_{qq}^3 - 10 DF_q F_{qq}^2 F_{qqq} + 10 F_p F_{qq}^2 F_{qqq} - 10 F_{qq}^4 F_z + 10 F_q F_{qq}^2 F_{qqq}F_z~]~\tilde{\om}^1\tilde{\om}^3 +\nonumber\\&30 F_{qq}^4~\tilde{\om}^1\tilde{\om}^4 + [~30 DF_q F_{qq}^3 - 30 F_p F_{qq}^3 - 30 F_q F_{qq}^3 F_z~]~\tilde{\om}^1\tilde{\om}^5+\nonumber\\&[~4 F_{qqq}^2 - 3 F_{qq} F_{qqqq}~]~(\tilde{\om}^2)^2-10 F_{qq}^2 F_{qqq}~\tilde{\om}^2\tilde{\om}^3 + 30 F_{qq}^3~\tilde{\om}^2\tilde{\om}^5 - 20 F_{qq}^4~(\tilde{\om}^3)^2,\nonumber 
\eeq
where the tilded omegas are defined by
\begin{eqnarray*}
&\tilde{\om}^1=\der y -p\der x\\
&\tilde{\om}^2=\der z-F\der x-F_q(\der p-q\der x)\\
&\tilde{\om}^3=\der p-q\der x\\
&\tilde{\om}^4=\der q\\
&\tilde{\om}^5=\der x.
\end{eqnarray*}
Despite of its ugliness this formula may be useful if one wants to
write down the Cartan invariant forms $(\theta^i,\Om_\mu)$ and the
scalar invariants $a_1,a_2,....$ directly in terms of the function
$F=F(x,y,y',y'',z)$ and its derivatives.\\

We can summarize the above considerations in the following theorem.
\bt~\\
The conformal class of
$(3,2)$-signature metrics $G_{(3,2)}$ which are naturally defined on
the $J$ space parametrized by $(x,y,p,q,z)$ encodes 
an invariant information about a class of forms 
\beq
&\om^1=\der z-F(x,y,p,q,z)\der x\nonumber\\
&\om^2=\der y-p\der x\nonumber\\
&\om^3=\der p-q\der x\nonumber\\
&\om^4=\der q\nonumber\\
&\om^5=\der x.\nonumber
\eeq
associated with a second order Monge equation 
$$z'=F(x,y,y',y'',z)$$
satisfying $F_{qq}\neq 0$.\\
Among all 5-dimensional $(3,2)$-signature metrics the metrics
$G_{(3,2)}$ are distinguished by the requirement that their
$\underline{\bf so}(4,3)$-valued Cartan normal conformal connection is
{\rm reducible} to a $\underline{\tilde{\bf g}_2}$-valued Cartan connection $\om_{\tilde{G}_2}$. 
\et

Interestingly the conformal metrics $G_{(3,2)}$ are very rarely
conformal to Einstein metrics. Even weaker curvature conditions, which
are necessary for a metric to be conformal to Einstein, such as
e.g. {\it conformal C-space conditions} (see Ref. \cite{rod} for the
definition), are not always satisfied by the metrics
$G_{(2,3)}$. However there are examples of the
second order Monge equations which correspond to the conformally
Einstein metrics $G_{(3,2)}$. Below, we present one of such examples.\\

\noindent
{\bf Example 6}\\
Consider a second order Monge equation
$$z'=F(y''),\quad\quad\quad\quad{\rm with}\quad\quad\quad\quad
F_{y''y''}\neq =0.$$
Since $F$ depends on only one variable $q$ we will denote its
derivatives by $F_q=F'$, etc. Its corresponding forms on $J$ are 
$\om^1=\der z-F(q)\der x$, $\om^2=\der y-p\der x$, $\om^3=\der p-q\der x$, 
$\om^4=\der q$, $\om^5=\der x$, with the tilded forms appearing in (\ref{32met})given by 
\begin{eqnarray}
&\tilde{\om}^1=\der y -p\der x\nonumber\\
&\tilde{\om}^2=\der z-F\der x-F'(\der p-q\der x)\nonumber\\
&\tilde{\om}^3=\der p-q\der x\label{przz}\\
&\tilde{\om}^4=\der q\nonumber\\
&\tilde{\om}^5=\der x.\nonumber
\end{eqnarray}
The invariant forms
$(\theta^1,\theta^2,\theta^3,\theta^4,\theta^5,\Om_1,\Om_2,\Om_3,\Om_4,\Om_5,
\Om_6,\Om_7,\Om_8,\Om_9)$ of theorem 8 are totally determined by forms 
$(\theta^1,\theta^2,\theta^3,\theta^4,\theta^5,\Om_1,\Om_2,\Om_3,\Om_4,\Om_5,
\Om_6,\Om_7,\Om_8,\Om_9)$ on $J$ which satisfy system 
(\ref{sycart})-(\ref{syspp}). Staring with (\ref{przz}) we find that
on $J$ these forms can be represented by 
\beq
&\theta^1=\tilde{\om}^1\nonumber\\
&\theta^2=\tilde{\om}^2\nonumber\\
&\theta^3=-(F'')^{\frac{1}{3}}\tilde{\om}^3\label{spec}\\
&\theta^4=(F'')^{-\frac{1}{3}}~[~\tilde{\om}^5-\frac{1}{3}F^{(3)}(F'')^{-1}\tilde{\om}^3+
\frac{1}{30}(-3F''F^{(4)}+4F^{(3)2})(F'')^{-3}\tilde{\om}^2~]\nonumber\\
&\theta^5=-(F'')^{\frac{2}{3}}\tilde{\om}^4\nonumber,
\eeq
\beq
&\Om_1=0\nonumber\\
&\Om_2=\frac{1}{90}[-45F''F^{(3)}F^{(4)}+40F^{(3)3}+9(F'')^2F^{(5)}](F'')^{-5}\theta^2+
\frac{1}{30}[-3F''F^{(4)}+4F^{(3)2}](F'')^{-\frac{10}{3}}\theta^3\nonumber\\
&\Om_3=0,\quad\quad\Om_4=0,\quad\quad\Om_5=0,\nonumber\\
&\Om_6=-\frac{1}{30}[-3F''F^{(4)}+4F^{(3)2}](F'')^{-\frac{10}{3}}\theta^5
\nonumber\\
&\Om_7=0,\quad\quad\Om_8=0,\quad\quad\Om_9=0.\nonumber
\eeq
In this setting the only nonvanishing function among
$(a_1,a_2,a_3,a_4,a_5,b_1,b_2,b_3,b_4,c_1,c_2,c_3,\delta_1,\delta_2,e,$
$h_1,h_2,h_3,h_4,$
$h_5,h_6,k_1,k_2,k_3)$ is 
\be
a_5=\{-224F^{(3)4}+336F''F^{(3)2}F^{(4)}-80(F'')^2F^{(3)}F^{(5)}+(F'')^2
[-51F^{(4)2}+10F''F^{(6)}]\}/[100(F'')^{\frac{20}{3}}].\label{a55}
\ee
Now applying formula (\ref{met32}) to the forms (\ref{spec}), or using
formula (\ref{32met}) for $F=F(q)$, we get the following representative
for the metrics $[G_{(3,2)}]$:
\beq
&G_{(3,2)}= 30(F'')^4~[~\der q\der y-p\der q\der
  x~]~+~[~4F^{(3)2}-3F''F^{(4)}~]~\der z^2+\nonumber\\&
2~[-5(F'')^2F^{(3)}-4F'F^{(3)2}+3F'F''F^{(4)}~]~\der p\der z+\nonumber\\&
2~[15(F'')^3+5q(F'')^2F^{(3)}-4FF^{(3)2}+4qF'F^{(3)2}+3FF''F^{(4)}-3qF'F''F^{(4)}~]~\der
  x\der z+\nonumber\\&
[-20(F'')^4+10F'(F'')^2F^{(3)}+4(F')^2F^{(3)2}-3(F')^2F''F^{(4)}~]~\der
  p^2+\label{metprzy}\\&
2~[-15F'(F'')^3+20q(F'')^4+5F(F'')^2F^{(3)}-10qF'(F'')^2F^{(3)}+\nonumber\\&
4FF'F^{(3)2}-4q(F')^2F^{(3)2}-3FF'F''F^{(4)}+3q(F')^2F''F^{(4)}~]~\der
  p\der x+\nonumber\\&
[-30F(F'')^3+30qF'(F'')^3-20q^2(F'')^4-10qF(F'')^2F^{(3)}+10q^2
  F'(F'')^2F^{(3)}+4F^2F^{(3)2}-\nonumber\\&
8qFF'F^{(3)2}+4q^2(F')^2F^{(3)2}-3F^2F''F^{(4)}+6qFF'F''F^{(4)}-3q^2(F')^2F''F^{(4)}~]~\der
x^2.\nonumber
\eeq
It is a matter of checking that this metric is conformal to an
Einstein metric $g={\rm e}^{2\Upsilon}G_{(3,2)}$ with the conformal
factor $\Upsilon=\Upsilon(q)$ satisfying equation 
$$10(F'')^2~[~\Upsilon''-(\Upsilon')^2~]-40F''F^{(3)}\Upsilon'+17F''F^{(4)}-56F^{(3)2}=0.$$

\noindent
Cartan \cite{5var} classified various types of nonequivalent forms
(\ref{cof}) according to the roots of the polynomial 
$$
\Psi(z)=a_1z^4+4a_2 z^3+6a_3 z^2+4 a_4 z+a_5,
$$
where $(a_1,a_2,a_3,a_4,a_5)$ are the scalar invariants given by
(\ref{syspp}). This polynomial encodes partial\footnote{For
  completeness we give the exact formula for The Weyl tensor of metrics
  $G{(3,2)}$ in the Appendix} information of the 
Weyl tensor of the associated metrics $G_{(3,2)}$. In particular, its
invariant $I_\Psi=6a_3^2-8a_2a_4+2a_1a_5$ is, modulo a numerical
factor, proportional to the square of the Weyl tensor
$C^2=C^{\mu\nu\rho\sigma}C_{\mu\nu\rho\sigma}$ of the metric
$G_{(3,2)}$. 
Vanishing of $I_\Psi$
means that $\Psi=\Psi(z)$ has a root with multiplicity no smaller than 3. Our
example above corresponds to the situation when this multiplicity is
equal to 4. According to Cartan \cite{5var} all nonequivalent forms
for which $\Psi$ has quartic root are covered by this example. The
nonequivalent classes are distinguished by the only nonvanishing
scalar invariant $a_5$ of (\ref{a55}), 
to which the Weyl tensor of metric (\ref{metprzy}) is proportional. 

\noindent
We were unable to construct an example of forms (\ref{cof}) for
which $\Psi$ has precisely triple root. For this it is enough to
assume that among the scalar invariants $(a_1,a_2,a_3,a_4,a_5)$ only
$a_4$ is nonvanishing. In such situation Cartan shows that the system
(\ref{sycart})-(\ref{syspp}) reduces to an invariant coframe on
$J$. Despite the fact that in this case the system is reducible to 5-dimensions it
is difficult, to find nonhomogeneous examples of forms which satisfy
it. 
\section{Acknowledgements}
This work was inspired by discussions on the
differential equations' aspect of the Null Surface Formulation of GR
which I had with Ezra Ted Newman. It started during my stay at King's College
London and was completed at the Institute for Mathematics of Humboldt
University in Berlin, where I was a member of the VW Junior
Research Group ``Special Geometries in Mathematical Physics''. 
I am very grateful to David C Robinson, 
Ilka Agricola and Thomas Friedrich for making the visits to these
institutions possible and fruitful.\\ 

\noindent
I acknowledge support from EPSRC grant
GR/S34304/01 while at King's College London and the KBN grant 2
P03B 12724 while at Warsaw University.
 
\section{Appendix}
In the null coframe $(\al^1,\al^2,\al^3,\al^4\al^5)=(\theta^1,\theta^2,\frac{2\sqrt{3}}{3}\theta^3,\theta^4,\theta^5)$ in which the metric
 (\ref{32met}) is 
$$G_{(3,2)}=2\al^1\al^5-2\al^2\al^4+(\al^3)^2$$
the Weyl tensor 2-forms are:
$$
C_{\mu\nu}=\frac{1}{2}C_{\mu\nu\rho\sigma}\al^\rho\dz\al^\sigma=
\bma 
0&0&0&-w_{14}&w_{15}\\
0&0&0&w_{15}&-w_{25}\\
0&0&0&-w_{34}&w_{35}\\
w_{14}&-w_{15}&w_{34}&0&-w_{45}\\
-w_{15}&w_{25}&-w_{35}&w_{45}&0
\ema,
$$
where 
$$
w_{14}=\frac{3}{8}c_3\al^1\dz\al^2+\frac{\sqrt{3}}{2}b_3\al^1\dz\al^3+a_3\al^1\dz\al^4+a_4\al^1\dz\al^5+\frac{\sqrt{3}}{2}b_4\al^2\dz\al^3+a_4\al^2\dz\al^4+a_5\al^2\dz\al^5,
$$
$$
w_{15}=\frac{3}{8}c_2\al^1\dz\al^2+\frac{\sqrt{3}}{2}b_2\al^1\dz\al^3+a_2\al^1\dz\al^4+a_3\al^1\dz\al^5+\frac{\sqrt{3}}{2}b_3\al^2\dz\al^3+a_3\al^2\dz\al^4+a_4\al^2\dz\al^5,
$$
$$
w_{25}=\frac{3}{8}c_1\al^1\dz\al^2+\frac{\sqrt{3}}{2}b_1\al^1\dz\al^3+a_1\al^1\dz\al^4+a_2\al^1\dz\al^5+\frac{\sqrt{3}}{2}b_2\al^2\dz\al^3+a_2\al^2\dz\al^4+a_3\al^2\dz\al^5,
$$
$$
w_{34}=\frac{3\sqrt{3}}{16}\delta_2\al^1\dz\al^2+\frac{3}{4}c_2\al^1\dz\al^3+\frac{\sqrt{3}}{2}b_2\al^1\dz\al^4+\frac{\sqrt{3}}{2}b_3\al^1\dz\al^5+\frac{3}{4}c_3\al^2\dz\al^3+\frac{\sqrt{3}}{2}b_3\al^2\dz\al^4+\frac{\sqrt{3}}{2}b_4\al^2\dz\al^5,
$$
$$
w_{35}=\frac{3\sqrt{3}}{16}\delta_1\al^1\dz\al^2+\frac{3}{4}c_1\al^1\dz\al^3+\frac{\sqrt{3}}{2}b_1\al^1\dz\al^4+\frac{\sqrt{3}}{2}b_2\al^1\dz\al^5+\frac{3}{4}c_2\al^2\dz\al^3+\frac{\sqrt{3}}{2}b_2\al^2\dz\al^4+\frac{\sqrt{3}}{2}b_3\al^2\dz\al^5,
$$
$$
w_{45}=-\frac{9}{64}e\al^1\dz\al^2+\frac{3\sqrt{3}}{16}\delta_1\al^1\dz\al^3+\frac{3}{8}c_1\al^1\dz\al^4+\frac{3}{8}c_2\al^1\dz\al^5+\frac{3\sqrt{3}}{16}\delta_2\al^2\dz\al^3+\frac{3}{8}c_2\al^2\dz\al^4+\frac{3}{8}c_3\al^2\dz\al^5.
$$

\end{document}